\documentclass[a4paper,11pt]{article}
\usepackage{amsmath,amssymb,mathtools}
\usepackage{parskip}
\usepackage{hyperref}

\newtheorem{thm}{Theorem}
\newtheorem{prop}[thm]{Proposition}
\newtheorem{lemma}[thm]{Lemma}
\newtheorem{cor}[thm]{Corollary}
\newtheorem{qn}{Question}
\newtheorem{conj}{Conjecture}

\begin{document}
\title{On profinite groups with positive rank gradient} 
\author{Nikolay Nikolov \footnote{Mathematical Institute, Oxford, OX2 6GG, UK. email: nikolov@maths.ox.ac.uk}} \date{MSC 20E18, 20F69}
\maketitle 

\begin{abstract}We prove that a profinite group $G$ with positive rank gradient does not satisfy a group law. In the case when $G$ is a pro-$p$ group we show that $G$ contains a nonabelian dense free subgroup. \end{abstract}
\section{Introduction}
Let $G$ be a finitely generated group. We denote by $d(G)$ the minimal number of generators of $G$, with the convention that if $G$ is a profinite group then $d(G)$ is the minimal number of topological generators. When $p$ is a prime $d_p(G)$ will denote the minimal number of generators of the pro-$p$ completion of $G$, that is $d_p(G)= \dim_{\mathbb F_p} \frac{G}{\Phi(G)}$, where $\Phi(G)=G^p[G,G]$ is the $p$-Frattini subgroup of $G$. 

Assume additionally that $G$ is residually finite. Let $G \geq H_1 \geq H_2 \geq \cdots$ be a chain of normal subgroups $H_i$, each having finite index in $G$ and such that $\cap_i H_i = \{1\}$. The rank gradient of $G$ with respect to $(H_i)$ is defined as \[ RG(G, (H_i)):=\lim_{i \rightarrow \infty} \frac{d(H_i)-1}{|G:H_i|}.\] 
When $(H_i)$ is a  chain of normal subgroups as above and additionally $|G:H_i|$ is a power of $p$ for all $i \in \mathbb N$, we define 
 the $p$-gradient $RG_p(G, (H_i))$ in the same way: $RG_p(G, (H_i)):=\lim_{i \rightarrow \infty} \frac{d_p(H_i)-1}{|G:H_i|}$.

From Schreier's inequality $d(H_i)-1 \leq |G:H_i| (d(G)-1)$ and its pro-$p$ analogue we deduce that $RG(G,(H_i))$ and $RG_p(G,(H_i))$ exist as a limits of monotonic non-increasing sequences.

 It is an open problem whether the rank gradient depends on the choice of the chain $(H_i)$ and this is related to the Fixed Price problem of topological dynamics, see \cite{AN}. In any case we set
\[ RG(G):= \inf  \left \{ \frac{d(H)-1}{|G:H|} \ | \ H<_f G  \right \},\] where the infimum is taken over all subgroups $H$ of finite index $G$, and define $RG_p(G)$ similarly:
\[ RG_p(G):= \inf  \left \{ \frac{d_p(H)-1}{|G:H|} \ | \ H \vartriangleleft  G, \ |G:H|=p^k, \ k \in \mathbb N   \right \}.\] 

Suppose $M$ is a subgroup of finite index of $G$. An easy exercise with the Schreier inequality gives that the infimum in the definition of $RG(G)$ can be restricted to the set of finite index subgroups $H$ of $G$ contained in $M$. Therefore we have $RG(M)= |G:M| RG(G)$.

When $G$ is a profinite group we define the rank gradient $RG(G,(H_i))$ of $G$ with respect to chains of open normal subgroups $H_i$ with $\cap H_i=\{1\}$ using the same expression as above.  In that case a compactness argument shows that $RG(G,(H_i))$ does not depend on the choice of the sequence $(H_i)$ and is equal to $RG(G)$.

Starting with the work of Lackenby \cite{L2}, there has been a lot of interest in rank gradient of abstract residually finite groups, e.g. \cite{AN}, \cite{AJN}, \cite{L}. It is tempting to believe that finitely generated abstract groups with positive rank gradient are in some way related to free groups. Notable progress in this was obtained by Lackenby \cite{L}, who proved that finitely presented residually $p$-groups with $RG_p(G)>0$ are large (meaning that a finite index subgroup has a free nonabelian homomorphic image). However Osin \cite{O} and Schlage-Puchta \cite{JC} constructed residually finite torsion groups with positive rank gradient, showing that the finite presentability condition in Lackenby's theorem cannot be omitted, and indeed the connection with free groups is not true in general.

In this paper we will focus on profinite groups, where the relationship between positive rank gradient and free groups is more compelling. Our main result is the following. 

\begin{thm} \label{main} Let $G$ be a finitely generated profinite group with positive rank gradient. Then $G$ does not satisfy a nontrivial group law. 
\end{thm}

A key ingredient in the proof is the following result of Dalla Volta and Lucchini \cite{DL}: the minimal size of a generating set of a finite group $G$ registers in some quotient $\bar G$ which is a \emph{crown-based power}, see Theorem \ref{crown} below.

Theorem \ref{main} is related to the following question of A. Thom.
\begin{qn}[A. Thom] \label{thom}
Let $\Gamma$ be a residually finite non-amenable group. Must the profinite completion $\hat \Gamma$ of $\Gamma$ contain a nonabelian free group? 
\end{qn}
A positive answer to Question \ref{thom}  implies Theorem \ref{main} since a dense finitely generated subgroup $\Gamma$ of $G$ must be non-amenable (see \cite{AJN}, Theorem 5 or \cite{AN}, Theorem 3). Therefore if $\hat \Gamma$ contains a free group then $\Gamma$ (and hence $G$) cannot satisfy a non-trivial  group law.

A positive answer to Question \ref{thom} will also imply the following.

\begin{conj} Let $\Gamma$ be a finitely generated residually finite group which satisfies a non-trivial group law. Then $RG(\Gamma)=0$.  \end{conj}

While we could not prove this conjecture, Theorem \ref{main} implies the weaker result.
\begin{cor} \label{residual} Let $\Gamma$ be a finitely generated residually finite group which satisfies a non-trivial group law. Then
\[ \inf \left \{ \frac{d(\Delta^{ab})-1}{|\Gamma: \Delta|} \  | \ \Delta<_f \Gamma \right \} =0,\]
where $\Gamma^{ab}=\Gamma/[\Gamma, \Gamma]$ denotes the abelianization of $\Gamma$ and the infimum is taken over all subgroups $\Delta$ of finite index in $\Gamma$.
\end{cor} 

Indeed, let $G=\hat \Gamma$ be the profinite completion of $\Gamma$. Since $G$ satisfies the group laws of $\Gamma$, Theorem \ref{main} gives $RG(G)=0$. In particular $\frac{d(G_i)-1}{|G : G_i|} \rightarrow 0 $ for some sequence of open subgroups $G_i$ of $G$. Let $\Gamma_i = \Gamma \cap G_i$ and observe that $G_i$ is the profinite completion of $\Gamma_i$, hence $d(\Gamma_i^{ab}) =d((\hat \Gamma_i)^{ab}) \leq d(G_i)$. Corollary \ref{residual} follows.

Question \ref{thom} has the following profinite version.  

\begin{qn} \label{free}
Must a profinite group with positive rank gradient contain a nonabelian free subgroup? 
\end{qn}

A \emph{coset identity} on a group $G$ is a nontrivial group word $w(x_1, \ldots, x_k)$ in $k$ letters, a normal subgroup $H$ of  finite index in $G$ and $k$ cosets $g_1H,\ldots ,g_kH$ from $G/H$, such that $w(g_1h_1, \ldots, g_kh_k)=1$ for all $h_1, \ldots, h_k \in H$. In case $G$ is a profinite group we also require that $H$ is open in $G$.

An affirmative answer to Question \ref{free} will follow if one can prove that a finitely generated profinite group with positive rank gradient does not satisfy a coset identity. 
While the author believes that this is true, the method of proof of Theorem \ref{main} does not give it. 
We will prove the following.

\begin{thm} \label{relative} Let $G$ be a finitely generated profinite group which satisfies a coset identity of length $m$. Assume that $d(G)>6m$. Then $G$ has an open normal subgroup $H$ with $\frac{d(H)-1}{|G:H|}<\alpha (d(G)-1)$, where $\alpha= \frac{2m-1}{2m}$. 
\end{thm}

Theorem \ref{main} is easily deduced by repeated applications of Theorem \ref{relative}. We cannot use the same method to answer Question \ref{free} because, unlike group laws, coset identities may not  induce coset identities on subgroups. \medskip

We prove an affirmative answer to Question \ref{free} in the case of pro-$p$ groups.

\begin{thm} \label{prop} Let $G$ be a finitely generated pro-$p$ group with positive rank gradient. Then $G$ contains a dense non-abelian free subgroup.
\end{thm}
The proof of Theorem \ref{prop} is much shorter than the proof of Theorem \ref{main}. It relies on an application of Schlage-Puchta's result \cite{JC}, together with Lie algebra methods originally developed by Wilson and Zelmanov \cite{WZ} to study Golod-Shafarevich groups, as well as Zelmanov's results \cite{Z} on Lie algebras with identities.

\section{Proofs}

All modules in the paper will be left modules and to be consistent we will write group actions on the left as well. This applies in particular to conjugation: for two elements $a,b$ of a group $G$ we write $\prescript{a}{} b:= aba^{-1}$. For sets $X$ and $Y$ we write $X \backslash Y$ for the elements of $X$ outside $Y$. By way of contrast $G/H:= \{ xH \ | \ x \in G\}$ will denote the left cosets of a subgroup $H$ of a group $G$.

\subsection{Proof of Theorem \ref{prop}}

Recall the following useful result of Levai and Pyber. We add a proof of it for completeness.

\begin{prop}[Levai-Pyber 1998]  \label{cosetid} A pro-$p$ group $G$ contains no dense non-abelian free subgroups if and only if $G$ satisfies a coset identity.
\end{prop}
\textbf{Proof:} A coset identity of $G$ induces a coset identity on its dense subgroups. Therefore if $G$ satisfies a coset identity $G$ has no dense non-abelian free subgroups. 

Conversely, suppose $G$ has no dense non-abelian free subgroups. Let $g_1, \ldots, g_k$ be a topological generating set of $G$; without loss of generality we can take $k \geq 2$. Let $Y= (g_1 \Phi(G)) \times \cdots \times (g_k\Phi(G))$. Then $Y$ is a closed and open subset of $G^k=G \times \cdots \times G$ ($k$ times). Note than any tuple $\mathbf y=(y_1 , \ldots ,y_k) \in Y$  generates a dense subgroup $\Gamma_\mathbf y:=\langle y_1, \ldots, y_k \rangle$ of $G$. Since $\Gamma_\mathbf y$ is not a free group there is a reduced word $w \not =1$ in the free group $F_k$ on $k$ letters, such that $w(y_1, \ldots, y_k)=1$.

It follows that $Y= \cup_{ w \in F_k \backslash \{1\}} Z_w$, where 
\[ Z_w= \{ \mathbf y=(y_1, \ldots, y_k) \in Y \ | \ w(y_1, \ldots y_k) =1\}.\]
Each $Z_w$ is a closed subset of $Y$. Since $Y$ is compact, Hausdorff and a union of countably many closed sets $Z_w$, the Baire category theorem implies that some $Z_{w_0}$ contains a nonempty open set. The base of the topology of $G$ is given by cosets of open normal subgroups. Therefore there is some open normal subgroup $H$ of $G$ and cosets $u_1H, \ldots, u_kH$ such that $(u_1H) \times \cdots \times (u_k H) \subseteq Z_{w_0}$. This means $w_0(u_1H, \ldots, u_kH)=1$ and so $G$ satisfies a coset identity.
$\square$

For a group $\Gamma$ and a prime $p$ we recall $(D_n(\Gamma))_{n=1}^\infty$, the Zassenhaus filtration of $\Gamma$, defined by $D_n(\Gamma)= \{g \in \Gamma \ | \ g-1 \in I^n \}$, where $I$ is the augmentation ideal of the group algebra $\mathbb F_p\Gamma$. The Lie algebra $\mathcal L_p(\Gamma)$ is defined as 
\[ \mathcal L_p(\Gamma)= \bigoplus_{i=0}^\infty \frac{D_i(\Gamma)}{D_{i+1}(\Gamma)}\] with Lie bracket $[gD_{i+1}(\Gamma),g'D_{j+1}(\Gamma)]= [g,g']D_{i+j+1}(\Gamma)$ for all $g \in D_i(\Gamma)$ and $g' \in D_j(\Gamma)$. 
The following was proved in \cite{WZ}.
\begin{thm} \label{liealg} Let $\Gamma$ be a residually finite $p$-group which satisfies a coset identity. Then $\mathcal L_p (\Gamma)$ satisfies a Lie algebra identity.
\end{thm}

We shall need the following variation of a result of Schlage-Puchta proved in \cite{JC} in the case when $\Gamma$ is a free group.
\begin{thm} \label{puchta} Let $p$ be a prime integer and let $\Gamma$ be a finitely generated group with an infinite chain $\Gamma>\Gamma_1> \cdots $ of normal subgroups of $p$-power index in $\Gamma$. Let \[ \hat \Gamma = \underleftarrow{\lim}(\Gamma/\Gamma_i)\] be the completion of $\Gamma$ with respect to $(\Gamma_i)$ and assume that $RG_p(\hat \Gamma)>0$ Then $\Gamma$ has a quotient $\Delta$ with a chain $(\Delta_i)$ such that $RG_p(\Delta, ( \Delta_i))>0$ and $\Delta$ is an infinite residually finite $p$-torsion group.  
\end{thm}

This theorem can be deduced easily from the results in \cite{JC} and for completeness we give a proof of it in the next subsection.  We continue with the proof of Theorem \ref{prop}.
Suppose that the claimed result is false and $G$ is a pro-$p$ group which is a counterexample. By Proposition \ref{cosetid} $G$ satisfies a coset identity. Let $(G_n)$ be a chain of normal subgroups with trivial intersection in $G$. We choose a finitely generated dense subgroup $\Gamma$ inside $G$ and let $\Gamma_i= G_i \cap \Gamma$.  For each $i\in \mathbb N$ we have $G/G_i \simeq \Gamma/\Gamma_i$ and therefore the profinite completion of $\Gamma$ with respect to $(\Gamma_i)$ is $G$.

Since $RG_p(G)>0$ Theorem \ref{puchta} implies that $\Gamma$ has an infinite $p$-torsion quotient $\Delta$.
On the other hand, since $G$ satisfies a coset identity then so does $\Gamma$ and its quotient $\Delta$. Therefore by Theorem \ref{liealg} $\mathcal L_p(\Delta)$ satisfies a Lie algebra identity. In addition, every homogeneous element of $\mathcal L_p(\Delta)$ is ad-nilpotent since $\Delta$ is $p$-torsion. By Theorem 1.1 of \cite{Z} the Lie algebra $\mathcal L_p(\Delta)$ is finite dimensional and hence $\Delta$ must be finite, a contradiction. Theorem \ref{prop} follows. $\square$

\subsection{Proof of Theorem \ref{puchta}}
Let $F$ be a nonabelian free group and $g \in F \backslash \{1\}$. We define $\nu_{F}(g)$ as the largest integer $k \geq 0$ such that $g=h^{p^k}$ for some $h \in F$. For completeness we set $\nu_F (1)= \infty$.

For a subset $X \subseteq F$ define \[ \delta_F(X) = \sum_{g \in X} p^{-\nu_F(g)}, \] with the convention that $\delta_F(X)= \infty$ if the above sum diverges.

Observe that if $g \in X$ with $\nu_F(g) \geq 1$ then $g \in F^p \leq \Phi(F)$, where $\Phi(F)=F^p[F,F]$ is the $p$-Frattini subgroup of $F$. In particular $p^{-\nu_F(g)}=1$ for every $g  \in X \backslash \Phi(F)$, giving 
\begin{equation} \label{dp} |X \backslash \Phi(F)| \leq \delta_F(X). \end{equation}

The following technical result is proved in \cite{JC} as part of the proof of Theorem 2 there.

\begin{lemma} \label{subgroup} Let $F$ be a nonabelian free group and let $N= \langle  \prescript{F}{} X \rangle \leq F$ be a subgroup which is the normal closure in $F$ of a set $X$. Let $H$ be a subnormal subgroup of $F$ with $H \geq N$ and $|F:H|=p^n$ for some $n \in \mathbb N$. Then $N$ contains a subset $Y \subseteq  \prescript{F}{} X$ such that $N=\langle \prescript{H}{} Y\rangle$ is the normal closure of $Y$ in $H$ and $\delta_H(Y) \leq |F:H| \delta_F(X)$.  
\end{lemma}

We are now ready to prove
Theorem \ref{puchta}. Adopt the notation and hypothesis of the theorem and let $\epsilon = RG_p(\hat \Gamma)$. Let $F$ be a free group projecting onto $\Gamma$. Without loss of generality we may assume that $\Gamma = F/M$, where $M$ is a normal subgroup of $F$. Let $F_i$ be the normal subgroups of $F$ such that $\Gamma_i=F_i/M$. We define the open subgroups of $F$ to be the preimages of the open subgroups of $\hat \Gamma$ under the composition $F \rightarrow \Gamma \rightarrow \hat \Gamma$. Specifically a subgroup $H$ of $F$ is open if $H$ contains some $F_i$.

Let $w_1, w_2, \ldots $ be an enumeration of the elements of $F \backslash \{1\}$ and define 
\[ X= \{ w_1^{p^{m_1}}, w_2^{p^{m_2}}, \ldots  \},\]
where $m_1<m_2< \ldots $ is a sequence of integers such that $\sum_{i=1}^\infty {p^{-m_i} }< \epsilon/2$.  

Let $N=\langle ^ F X\rangle$ be the normal closure of $X$ in $F$ and note that by construction we have $\delta_F (X) < \epsilon/2$ and $F/N$ is a $p$-torsion group.

Let $\tilde \Gamma:=F/MN$ and define the open subgroups of $\tilde \Gamma$ to be the images of the open subgroups of $F$, i.e. a subgroup $\tilde H =\frac{H}{MN} \leq \tilde \Gamma$ is open in $\tilde \Gamma$ if $H$ contains some $F_i$. We will prove the following
\begin{prop} \label{boundd}Let $\tilde H$ be an open subgroup of $\tilde \Gamma$. Then there is an open subgroup $\tilde A$ of $\tilde \Gamma$ with $\tilde H > \tilde A \geq \Phi (\tilde H)$ and 
\[ d_p(\tilde H/ \tilde A) -1> \frac{1}{2} \epsilon |\tilde \Gamma : \tilde H |\]
\end{prop}
\textbf{Proof:}
Let $n= |\tilde \Gamma : \tilde H|$ and let $\tilde H= H/NM$ where $H \geq NM$ is an open subgroup of $F$ with $|F:H|=n$. Denote by $L$ the closure of  $H/M \leq \Gamma $ inside $\hat \Gamma$, this is an open subgroup of $\hat \Gamma$ of index $n$. Since $RG_p(\Gamma) =\epsilon$ we have $\frac{d_p(L)-1}{n} \geq  \epsilon$. 

We have $d_p(L)=\dim_{\mathbb F_p} (L/\Phi (L))$, where the Frattini subgroup $\Phi(L)$ is an open subgroup of $L$ and of $\hat \Gamma$.

Let $B/M = \Phi(L) \cap \Gamma$. Thus $B$ is an open subgroup of $F$ with $H \geq B \geq \Phi(H)M$ and $d_p(H/B)=\dim_{\mathbb F_p} (L/\Phi(L)) \geq 1+ \epsilon n$. 

By Lemma \ref{subgroup} applied to $N=\langle  ^F X \rangle$ and the open subgroup $H\geq N$ we deduce that $N=\langle  ^H Y \rangle$ is the normal closure in $H$ of a set $Y$ with $\delta_H(Y) \leq |F:H| \delta_F(X) < \epsilon n /2$.  Let $A=NB$, since $A \geq B$ this is an open subgroup of $F$ and  so $\tilde A= A/MN$ is an open subgroup of $\tilde \Gamma$.  We have $ H \geq A \geq \Phi(H) MN$ and hence $\tilde H \geq \tilde A \geq \Phi(\tilde H)$.  Let $Y_0$ be the subset of $Y$ outside $\Phi(H)$. By (\ref{dp}) we have $|Y_0| \leq \delta_H(Y) < \epsilon n/2$. Since $H/B$ is an elementary abelian $p$-group we have 
$A=NB=\langle Y_0 \rangle  B$ and
\[ d_p (\tilde H / \tilde A) = \dim_{\mathbb F_p} \frac{H}{A} \geq \dim_{\mathbb F_p} \frac{H}{B} -|Y_0| > 1+ \epsilon n - \frac{\epsilon n}{2} = 1+ \frac{\epsilon n}{2}.\]
Observe that since $d_p(\tilde H/\tilde A)>1$ it follows that $\tilde H > \tilde A$. Proposition \ref{boundd} follows.
$\square$ \medskip

We can apply Proposition \ref{boundd} successively to find a sequence $\tilde \Gamma_1=\tilde \Gamma > \tilde \Gamma_2 > \cdots $ of open subgroups of $\tilde \Gamma$ such that $\tilde \Gamma_i > \tilde \Gamma_{i+1} \geq \Phi(\tilde \Gamma_i)$ and \[ \frac{d_p(\tilde \Gamma_i/\tilde \Gamma_{i+1}) -1}{|\tilde \Gamma : \tilde \Gamma_i|}> \epsilon/2 \] for all $i \in \mathbb N$.
Noting that $\tilde \Gamma$ is a $p$-torsion group Theorem \ref{puchta} follows by setting $\Lambda= \tilde \Gamma /U$ and  
$\Lambda_i= \tilde \Gamma_i /U$ where $U= \cap_{i=1}^\infty \tilde \Gamma_i$.  $\square$ 

\subsection{Proof of Theorem \ref{main}}

 Suppose $G$ satisfies a non-trivial group law of length $m$. We will show that $RG(G)=0$. This is clear if $G$ has a subsequence of subgroups $(G_i)_{i=1}^\infty$ such that $d(G_i)$ is bounded and $|G:G_i| \rightarrow \infty$. Therefore we may assume that $d(H) \rightarrow \infty$ as $H$ ranges over all open subgroups of $G$ and $|G:H| \rightarrow \infty$. In particular there is an open subgroup $G_0 \leq G$ such that $d(H)>6m$ for all open subgroups $H \leq G_0$. Since $RG(G_0)=|G:G_0|RG(G)$ it is sufficient to show $RG(G_0)=0$. Hence by replacing $G$ with $G_0$ we may assume that $d(H)>6m$ for any open subgroup $H$ of $G$. 

Since a group law is also a coset identity of $G$, Theorem \ref{relative} gives that $G$ has an open normal subgroup $H$ with $\frac{d(H)-1}{|G:H|} < \alpha (d(G)-1)$ where $\alpha= \frac{2m-1}{2m}$.

By replacing $G$ with $H$ and iterating Theorem \ref{relative}  $n$ times, we obtain a sequence of open subgroups $H_n<H_{n-1} < \cdots H_1<H_0=G$ such that $\frac{d(H_i)-1}{|H_{i-1}:H_i|}< \alpha (d(H_i)-1)$ for each $i$, giving that $\frac{d(H_n)-1}{|G:H_n|}< \alpha^n d(G)$. Since $\alpha<1$ by letting $n \rightarrow \infty$ we obtain $RG(G) \leq RG(G,(H_i))=0$.

\subsection{Proof of Theorem \ref{relative}}

Without loss of generality we may assume that the length $m$ of the coset identity on $G$ is minimal possible. 
Let $F(X)$ be the free group on $X=\{x_1,x_2, \ldots, x_k\}$ and let $w(x_1, \ldots, x_k) \in F(X)$ be a reduced group word of minimal length $m$ which is a coset identity on $G$. Thus there is an open normal subgroup $H$ of $G$ and $g_1, \ldots, g_k \in G$ such that $w(g_1H, \ldots, g_kH)=1$.

Let $w(x_1, \ldots, x_k)=w_1w_2  \cdots w_{m}$, where  $w_i=x_{t_i}^{\epsilon_i}$ for $i=1, \ldots, m$ with $\epsilon_i \in \{\pm 1\}$ and $x_{t_i} \in X$. Put $u_j=w_1 \ldots, w_j$ for $j=1, \ldots, m$. We will refer to $u_j$ as the initial subwords of $w$. Without loss of generality we may assume $t_1=\epsilon_1=1$ so that $u_1=w_1=x_1$. Moreover $w_{m} \not= w_1^{-1}=x_1^{-1}$, otherwise $w$ will be conjugate to the shorter identity $w_2 \cdots w_{m-1}$.
  
  We claim that by passing to a $G$-normal subgroup of $H$  and replacing $g_i$ with appropriate elements from $g_jH$ we may assume that the evaluations of the initial subwords $h_j:=u_j(g_1, \ldots g_k)$ are non-trivial for $j=1,2, \cdots m-1$. 

Indeed, suppose that we have obtained that $h_1, \ldots, h_{r-1} \not =1$ for some $r<m$.
Replace $H$ with a smaller open normal subgroup of $G$ such that $h_j \not \in H$ for $j=1, \ldots, r-1$. For any choice of $x_i \in H$ ($i=1, \ldots, k$) we have $u_j(g_1 x_1, \ldots, g_kx_k) \equiv  h_j$ mod  $H$ and in particular $u_j(g_1 x_1, \ldots, g_kx_k) \not =1$ for $j=1, \ldots, r-1$. As $r <m$, by the minimality of $m$ the word $u_{r}$ is  not a coset identity of $G$ and therefore we can find elements $x_i \in H$ with $u_{r}(g_1 x_1, \ldots, g_kx_k) \not =1$. Replace $g_i$ with $g_ix_i$ and we have achieved $h_{j} \not = 1$ for $j=1, \ldots, r$.

We can repeat this procedure increasing $r$ one at a time until we obtain that $h_1, \ldots, h_{m-1} \not =1$, proving the claim. \medskip

Next, by replacing $H$ with a smaller open normal subgroup of $G$ we may assume that $h_j \not \in H$ for $j=1,2, \ldots, m-1$.

Choose an element $y \in H$ and write $w(yg_1,g_2, \ldots, g_k)= u \cdot w(g_1, \ldots, g_k)$ where $u =u(y) \in H$ can be written as a product of several conjugates of $y^{\pm 1}$:

\[ u(y) = y \cdot \prescript{b_1}{} y^{e_1} \cdots  \prescript{b_s}{} y^{e_s}. \]
Here $s \geq 0$ is the number of occurences of $x_1^{\pm 1}$ among the $w_2, \ldots, w_m$, $e_i \in \{\pm 1\}$ and $(b_1, \ldots, b_s)$ is a subsequence of $(h_1, \ldots, h_{m-1})$.
Notice that we use the fact that $w_m \not = x_1^{-1}$ here to remove the possibility that the last $b_s$ could end up equal to $h_m=w(g_1, \ldots, g_k)=1$.
Since $H$ is normal in $G$ and $y \in H$, we have $yg_1 \in g_1H$ and therefore $w(yg_1,g_2, \ldots, g_k)=w(g_1,g_2, \ldots, g_k)=1$. We conclude that $u(y)=1$ for all $y \in H$.

We now arrive at the main technical result of this paper:

\begin{thm} \label{iden} Let $G$ be a profinite group with an open normal subgroup $H$ and elements $b_1, \ldots, b_s \in G \backslash H$. Let $e_i \in \{\pm 1\}$ for $i=1, \ldots, s$. For any $y \in H$ define $u(y)= y \cdot \prescript{b_1}{} y^{e_1} \cdots \prescript{b_s}{} y^{e_s}$ and assume $u(y)=1$ for all $y \in H$. Assume that $d(G)> 6s+6$. 
Then \[ \frac{d(H)-1}{|G:H|} < \alpha (d(G)-1),\] where $\alpha=\alpha(s)=\frac{2s+1}{2s+2}$.
\end{thm}

Since $s \leq m-1$ we have $\alpha(s)  \leq \frac{2m-1}{2m}$ and Theorem \ref{relative} follows from Theorem \ref{iden} and the above discussion.

Note that if $s=0$  then $u(y)=y$ which implies that $H=\{1\}$ and the claimed inequality holds trivially. Therefore without loss of generality we may assume that $s \geq 1$ (and hence $H \not = G$) in the rest of the argument.

For a profinite group $H$ the minimal number of generators $d(H)$ is the minimum of $d(\bar H)$ for all finite topological images $\bar H$ of $H$. Hence in order to prove Theorem \ref{iden} we may assume that $H$ and $G$ are finite groups. We shall also need more information on the minimal quotients of $H$ which realize $d(H)$ and the notion of a \emph{crown-based power}.
 
\subsection{Crown-based powers}

Let $L$ be a non-cyclic finite group with a unique minimal normal subgroup $N$. If $N$ is abelian then we also require that $N$ has a complement in $L$. For $k \geq 2$ the crown-based power $L_k$ of $L$ is defined as 

\[ L_k=\{ (a_1,\ldots, a_k) \in L^k \ | \ a_1N=a_2N= \cdots a_kN \},\]
where $L^k=L \times L \times \cdots \times L $ ($k$ times). For later use we define the projections
\begin{equation} \label{ti} t_i: L_k \rightarrow L, \quad t_i(a_1, \ldots, a_k)=a_i, \quad  \forall (a_1, \ldots, a_k ) \in L_k , \ i =1, \ldots, k \end{equation}
and note that $t_1(x) \equiv \cdots  \equiv t_k(x)$ mod $N$ for every $x \in L_k$.

The precise relationship between $d(L_k)$ and the integer $k$ is established in the following theorem from \cite{VL}.
When $N$ is nonabelian by  $P_{L,N}(d)$ we denote the probability that $d$ random elements of $L$ generate it, provided they generate $L/N$. By Theorem 1.1 of \cite{DL} we have $53/90 \leq P_{L,N}(d)\leq 1$.

When $N$ is abelian then $N$ is a simple faithful $L/N$-module and we denote $E= End_{L/N}(N)$, a finite field.
\begin{thm}[\cite{VL}, Theorem 2.7] \label{gen} Let $L_k$ be a crown-based power of $L$ and choose an integer $d \geq d(L)$. Then $d(L_k) \leq d$ if and only if one of the following holds: \medskip

1. The group $N$ is abelian and $k\leq (d-1) \dim_E N - \dim_E H^1(L/N,N)$.
\medskip

2. The group $N$ is nonabelian and $k\leq P_{L,N}(d) |N|^d /|C_{Aut(N)}(L/N)|$.
\end{thm}

Crown-based powers are the minimal among the $d$-generated finite groups in the following sense.
 
\begin{thm}[\cite{VL}, Theorem 1.4] \label{crown} Let $H$ be a finite group  with $d(H)=d > 2$. Then $H$ has a normal subgroup $M<H$ such that $H/M$ is isomorphic to some crown-based power $L_k$ as above with $d(L_k)=d>d(L)$. 
\end{thm}

Before we start the proof of Theorem \ref{iden} we recall some standard notation. Let $A$ be a group with a (left) action of a group $H$ on $A$ by automorphisms. By a $H$-chief factor of $A$ we mean a section $Y/X$, where $X$ and $Y$ are normal $H$-invariant subgroups of $A$ such that there are no $H$-invariant normal subgroups $Z$ of $A$ with $X< Z <Y$. Further, we say that two $H$-groups $A$ and $B$ are $H$-isomorphic if there is a group isomorphism $\pi : A \rightarrow B$ which is compatible with the action of $H$, i.e. $h \cdot \pi (a)= \pi ( h \cdot a)$ for all $h \in H$ and $a \in A$. 

When $A=H$ and the action of $H$ is given by conjugation we shall refer to the $H$-chief factors of $H$ as the  \emph{chief factors} of $H$. We will need two additional technical results.

\begin{prop}\label{basic} Let $H$ be a finite group with a normal subgroup $U$. Let $K_i \leq U$ be $k$ distinct normal subgroups of $H$ such that $U/K_i$ is a chief factor of $H$ for $i=1, \ldots, k$.
Then  C1 imples C2 implies C3 below.

C1. Each $U/K_i$ is non-abelian.

C2. For each $i=1,2 , \ldots, k$ we have $K_i \not \geq  \cap_{j \not = i} K_j$.

C3. The diagonal homomorphism $U \rightarrow \prod_{i=1}^k (\frac{U}{K_i})$ induces an isomorphism \[ \frac{U}{\cap_i K_i} \rightarrow \prod_{i=1}^k \frac{U}{K_i} .\]
\end{prop}

\textbf{Proof:} We will first show that C2 implies C3. Let $R$ be the intersection of $K_2, \ldots, K_k$. By assumption $K_1 \not \geq R$ and hence $K_1<K_1R=U$ by the maximality of $K_1$ in $U$. Therefore \[ \frac{U}{ \cap_{i=1}^k K_i}= \frac{U}{K_1 \cap R} \simeq \frac{U}{K_1} \times \frac{U}{R}.\] The subset $K_2, \ldots, K_k$ satisfies C2 as well and  hence the conclusion follows by induction on $k$.

Next we show that C1 implies C2. Indeed, suppose that $K_1$ contains $\cap_{i=2}^k K_i$ and choose a subset $X \subseteq \{K_2, \ldots, K_{k}\}$ minimal with respect to $K_1 \geq \cap_{K_i \in X} K_i$. Then $X$ satisfies C2 and therefore the diagonal homomorphism induces an isomorphism \[ \frac{U}{\bigcap_{K_i \in X} K_i} \rightarrow \prod_{K_i \in X}\frac{U}{K_i}.\]
   
Let $\bar K$ be the image of $K_1/(\cap_{K_i \in X} K_i)$ in $ \prod_{K_i \in X}(U/K_i)$. Since $K_1K_i=U$ for each $K_i \in X$ we have that $\bar K$ projects onto each direct factor $U/K_i$ of  $\prod_{K_i \in X}(U/K_i)$. In particular since $\bar K$ is a normal subgroup it follows that \[ \bar K \geq [\bar K, (U/K_i)]=[(U/K_i),(U/K_i)]=U/K_i \] for each $K_i \in X$, since $U/K_i$ is perfect being a nonabelian chief factor of $H$.  

Therefore $\bar K=  \prod_{K_i \in X}(U/K_i)$ giving $K_1=U$, contradiction. Hence C1 implies C2 as claimed.
$\square$ \medskip

We will also need the following result which detects crown-based powers.

\begin{prop} \label{spot} Let $L$ be a noncyclic finite group with a unique minimal normal subgroup $N$. If $N$ is abelian assume additionally that $N$ is complemented in $L$.
Let $H$ be a finite group with a normal subgroup $U$. Let $K_i \leq U$ for $i=1, \ldots, k$ be pairwise distinct normal subgroups of $H$ together with an isomorphism $\alpha: H/K_1 \rightarrow L$ such that $\alpha(U/K_1)=N$. 

Assume that for each $i = 1, 2 \ldots, k$ there is an isomorphism $\beta_i: H/K_i \rightarrow H/K_1$ such that 
$\beta_i(hU/K_i)= h U/K_1$ for all $h \in H$. 

Then $H/ (\cap_{i=1}^k K_i)$ is isomorphic to a crown-based power $L_{k'}$ for some $k' \leq k$. Moreover, if $N$ is non-abelian then $k'=k$.
\end{prop}

\textbf{Proof:} 
If there is some $i \in \{1, \ldots, k\}$ such that $K_i \geq \cap_{j \not=i} K_j$ then we can omit $K_i$ from our list of normal subgroups and the statement to be proved remains unaffected. Therefore we may assume that each $K_i$ does not contain the intersection of the rest of the $K_j$.
By Proposition \ref{basic} this is automatically the case if $N$ is non-abelian. 

For each $i=1, \ldots, k$ let $\tilde \beta_i: H \rightarrow L$ be the composition $\tilde \beta_i (h)= \alpha \circ \beta_i (hK_i)$ with kernel $K_i$. and we note that \begin{equation} \label{u} \tilde \beta_1 (h)N= \tilde \beta_2 (h)N= \cdots = \tilde \beta_k(h)N  \quad \forall h \in H \end{equation}
Consider the map $\delta: H \rightarrow L^k$ defined by $\delta(h)=(\tilde \beta_1(h), \ldots, \tilde \beta_k(h))$ for all $h \in H$ and note that $\ker \delta= \cap_i K_i$. Moreover the image of $\delta$ is a subgroup of $L_k$ by (\ref{u}) above.

By Proposition \ref{basic} we have $U/(\cap K_i) \simeq \prod_i (U/K_i)$ since condition C2 holds.  Comparing sizes of $H/ (\cap_i K_i)$ and $L_k$ we deduce that $\delta(H)= L_k$ as required.

$\square$
 
\subsection{Proof of Theorem \ref{iden}: Reduction}

Suppose first that $G$ has a subgroup $N$ with $N \geq H$ and $d(N) < d(G)$. In particular $|G:N| \geq 2$.
Using the Schreier inequality and the fact that $\alpha(s)>\frac{1}{2}$ we obtain
\[ \frac{d(H)-1}{|G:H|} \leq \frac{d(N)-1}{|G:N|}< \frac{d(G)-1}{2}< \alpha(s) (d(G)-1). \] Theorem \ref{iden} follows. 

From now on we shall assume that 
\begin{equation} \label{bd} d(N) \geq d(G) > 6s+6 \end{equation}
for any subgroup of $G$ with $N \geq H$.

In particular $d(H) \geq d(G)>2$ and hence Theorem \ref{crown} applies to $H$.
Let $M$ be the normal subgroup of $H$ provided by Theorem \ref{crown}. So $H/M$ is isomorphic to some crown-based power $L_k$ as above with $d(L_k)=d(H)$.

In this section we will reduce the Theorem to proving Proposition \ref{iden1} which assumes $M$ is normal in $G$. We will refer to the chief factors of $H$ as $H$-chief factors. 

Let $U/M$ be the socle (i.e. the product of the minimal normal subgroups) of $H/M$. Thus $U$ is normal in $H$ and $U/M$ is isomorphic to $N^k$ where each summand $N$ is an  $H$-chief factor. 

Let $J=N_G(U)$ and choose representatives $a_1=1, a_2, \ldots, a_r$ in $G$ for the left cosets $G/J$ of $J$.
Define $U_1=U, U_2= \prescript{a_2}{}U, \ldots, U_r= \prescript{a_r}{}U$, these are the distinct conjugates of $U$ in $G$.

We define \[ T_i = \bigcap _{g \in a_iJ} \prescript{g}{} M \quad  (i=1, 2,\ldots, r).\] In particular $T_1= \cap_{g \in J} \prescript{g}{} M$ is normalized by $J$ and $T_i= \prescript{a_i}{} T_1$. Proposition \ref{decompose} below will show that $T_1, \ldots, T_r$ are all distinct and in particular the action of $G$ by conjugation on $U_1, \ldots, U_r$ is the same as the action of $G$ on $T_1, \ldots, T_r$.

By definition $U_i/T_i$ is a subdirect product of all the $U_i/\prescript{g}{}M$ with $g \in a_iJ$
Moreover for $g \in a_iJ$ we have $H/ \prescript{g}{}M \simeq L_k$ with $U_i/\prescript{g}{}M \simeq N^k$. In particular $U_i / \prescript{g}{} M$ is a direct product of $k$ $H$-chief factors, each being a $g$-conjugate of a $H$-chief factor appearing in $U/M$. 

In particular when $N$ is abelian then $U_i/T_i$ is a semisimple $H/U_i$ module, all of whose simple factors are faithful $H/U_i$-modules. 

When $N$ is nonabelian observe that an element $l \in L \backslash N$ cannot act as an inner automorphism on $N$, otherwise we get that $C_L(N) \not =\{1\}$ contradicting the uniqueness of the minimal normal subgroup $N$ of $L$. 

We summarise the above discussion in

\begin{prop} \label{help} The $H$-chief factors of $U_i/T_i$ are contained in the union of the $H$-chief factors of  $U_i/ \prescript{g}{} M$ for $g \in a_iJ$. If $C$ is a $H$-chief factor of $U_i/T_i$ and $h \in H$ with $h \not \in U_i$ then conjugation by $h$ induces a non-inner automorphism of $C$. 
\end{prop}

We define $\mathfrak U = \cap_{i=1}^r U_i$ and $\mathfrak T = \cap_{i=1}^r T_i$.
Then $\mathfrak U$ and $\mathfrak T$ are normal subgroups of $G$ being the intersection of all $G$-conjugates of $U$ and $M$ respectively. \medskip

Consider the diagonal homomorphism $\mathfrak f: \mathfrak U \rightarrow \prod_{i=1}^r (U_i/T_i)$ defined by $\mathfrak f(y)=(yT_1, \ldots, y T_r)$ for each $y \in \mathfrak U$. We have $\ker \mathfrak f= \cap_{i=1}^r T_i=\mathfrak T$.

\begin{prop} \label{decompose} The homomorphism $\mathfrak f$ is surjective and induces an isomorphism $\mathfrak U/\mathfrak T \rightarrow \prod_{i=1}^r (U_i/T_i)$ with the conjugation action of $G$ on $\mathfrak U/\mathfrak T$ permuting the direct factors $U_i/T_i$ in the same way as $G$ acts on $U_1, \ldots, U_r$.
\end{prop} 
\textbf{Proof:} The result is trivial if $r=1$ and so we may assume $r>1$. In particular, since $U_i \not = U_j$ for $i \not = j$ and $|H:U_i|=|H:U_j|$ it follows that $\mathfrak U$ is a proper subgroup of $ U_i$ for each $i=1, 2 \ldots, r$.

We \textbf{claim} that $\mathfrak U T_i= U_i$ for each $i=1,2, \ldots, r$. 

Assume this for the moment. Therefore  $\mathrm{Im}(\mathfrak f)$  is isomorphic as an $H$-group to a subdirect product of all $U_i/T_i$ for $i=1, \ldots, r$. Observe that for $i \not =j$ every $H$-chief factor of $U_i/T_i$ is not $H$-isomorphic to any $H$-chief factor of  $U_j/T_j$. Indeed, whenever $C$ is a $H$-chief factor of $U_i/T_i$ Proposition \ref{help} gives that the kernel of the composition $H \rightarrow Aut(C) \rightarrow Out(C)$ is precisely $U_i$, and the groups $U_j$ and $U_j$ are distinct.
This proves that $U_1/T_1, \cdots, U_r/T_r$ do not share common $H$-chief factors. Now the Jordan-H\"older theorem gives that a subdirect product of these is in fact the full direct product. Therefore $\mathfrak f$ is surjective. The preimage of the direct factor $U_i/T_i$ under $\mathfrak f$ is $\mathfrak U \cap (\cap_{j \not =i}T_j)$. In particular the groups $T_1, \ldots, T_r$ are pairwise distinct. Thus the conjugation action of $G$ on $U_1, \ldots, U_r$ is the same as the conjugation action of $G$ on $T_1, \ldots, T_r$, in turn this is the same as the  conjugation action of $G$ on the $r-1$ element subsets of $\{T_1, \ldots, T_r\}$ and this results in the same action of $G$ on the direct factors of $\prod_{i=1}^r (U_i/T_i)$.

Proposition \ref{decompose} follows once we prove the claim above.

\textbf{Proof of claim:}
Suppose $\mathfrak U T_i < U_i$ and choose a subset $X \subseteq \{U_1, \ldots, U_r \}$ of minimal size  subject to $\cap_{U_j \in X} U_j=\mathfrak U$. Since $G$ acts transitively by conjugation on the $U_1, \ldots, U_r$ we may assume that $U_i \in X$. Since $r>1$ and $\mathfrak U \not = U_i$ we must have $|X|>1$. Let $R$ be the intersection of all members of $X \backslash \{U_i\}$, thus $R$ is a normal subgroup of $H$ such that $\mathfrak U=U_i \cap R$. By the minimality of $X$ we have $U_i \not \supseteq R$ and we choose  $h \in R \backslash U_i$. Now $[U_i,R] \leq  U_i \cap R=\mathfrak U$ and therefore $R$ acts trivially by conjugation on $Q=U_i/\mathfrak U T_i$. In particular $h$ centralizes every $H$-chief factor of $Q$ and those are a subset of the $H$-chief factors in $U_i/T_i$. This contradicts Proposition \ref{help} since we chose $h \not \in U_i$. Hence we must have $U_i= \mathfrak U T_i$ and the claim is proved. $\square$

Returning to the proof of Theorem \ref{iden} choose $y_0 \in U_1$ and let  $y \in \mathfrak U$ be such that $\mathfrak f(y)=(y_0T_1, 1,1, \cdots ,1) \in \prod_{i=1}^r (U_i/T_i)$. 

Recall that conjugation by an element $g \in G$ sends the  direct factor $U_1/T_1$ to a different direct factor $U_i/T_i$ of $\prod_{i=1}^r (U_i/T_i)$, unless $g \in N_G(U)=J$. Let $b_{i_1}, b_{i_2}, \ldots, b_{i_t}$ be the subsequence of those elements from $b_1, \ldots, b_s$ which lie in $J$, in the same order.
 
It follows that $\mathfrak f(u(y))=(u_0(y_0)T_1, *, \cdots, *)$, where we write $*$ for coordinates we are not interested in, and
\[ u_0(y_0):= y_0 \cdot \prescript{b_{i_1}}{} y_0^{e_{i_1}} \cdots \prescript{b_{i_{t}}}{}y_0^{e_{i_t}}.\] 

Since $u(y)=1$  we must have $u_0(y_0) \in T_1$ for all $y_0 \in U_1$.

Suppose that we show that $\frac{d(H) -1}{|J:H|} \leq \alpha (t) (d(J)-1)$. Then 

\begin{equation} \label{J} \frac{d(H)-1}{|G:H|} \leq \frac{ \alpha(t) (d(J)-1)|J:H|}{|G:H|}= \alpha(t)\frac{d(J)-1}{|G:J|} \leq \alpha(t)(d(G)-1),\end{equation}
where the last inequality follows from Schreier's bound $\frac{d(J)-1}{|G:J|} \leq d(G)-1$.
Theorem \ref{iden} then follows from (\ref{J}) since $t \leq s$ and $\alpha(t) \leq \alpha(s)$.

By (\ref{bd}) we may assume that $d(J)>6s+6$.

Therefore by replacing $G$ with $J$ and $u(y)$ with $u_0(y)$ we have reduced Theorem \ref{iden} to the following.
\begin{prop}\label{1/2} Let $G$ be a finite group with $d(G)>6s+6$ for some integer $s \geq 1$. Let $H$ be a normal subgroup of $G$, let $b_i \in G \backslash H$ and $e_i \in \{\pm 1\}$ for $i=1, \ldots, s$. Let $M <H$ be a normal subgroup of $H$ such that $H/M$ is isomorphic to a crown based power $L_k$ with $d(L_k)=d(H)>d(L)$. Let $U/M= \mathrm{soc}(H/M)$ and assume that $U$ is a normal subgroup of $G$.  Let $T= \cap_{g \in G} \prescript{g}{}M$ and assume that $u(y) \in T$ for all $y \in U$, where $u(y)= y \cdot \prescript{b_1}{} y^{e_1} \ \cdots \prescript{b_s}{} y^{e_s}$. 
Then \[ \frac{d(H)-1}{|G:H|} < \alpha(s) (d(G)-1),\]
where $\alpha(s)=\frac{2s+1}{2s+2}$.
\end{prop}
\medskip

From now on assume that $G$, $H$, $M$, $T$ and $u$ are as in the above Proposition.
Our next aim is to reduce the proof to the case where $M$ is normal in $G$, namely Proposition \ref{iden1} below. It will turn out that $U/T$ is a direct product of crown based powers of $L$ but these may be larger than $H/M \simeq L_k$. More precisely we have the following.
 
\begin{prop} \label{spl}
There is a normal subgroup $R$ of $H$ with $R \leq M$ such that 

1. We have $H/R \simeq L_{k'}$ with $k' \geq k$.

2. If $R=R_1,R_2 , \ldots, R_l$ are the distinct conjugates of $R$ in $G$ then 

\[ T= R_1 \cap \cdots \cap R_l \] and 
the diagonal homomorphism $f: U \rightarrow \prod_{i=1}^l (U/R_i)$ induces an isomorphism
\begin{equation} \label{fact} \frac{U}{T} \rightarrow \prod_{i=1}^l (U/R_i). \end{equation}
Under this isomorphism the conjugation action of $G$ on $U/T$ permutes the direct factors $U/R_i$ of $\prod_{i=1}^l (U/R_i)$ in the same way as $G$ permutes $R_1, \ldots, R_l$.
\end{prop}

\textbf{Proof:} We fix a surjection $\pi: H \rightarrow L_k$ with $\pi(U) = N^k$ and $\ker \pi =M$.
 
For $i=1, \ldots, k$ let $K_i = \ker (t_i \circ \pi)$, where $t_i: L_k \rightarrow L$ is the projection defined in (\ref{ti}).  It follows that each $K_i$ is a normal subgroup of $H$ such that $\cap_{i=1}^k K_i = M$.
From the definition of $L_k$ we have that $t_i$ is surjective onto $L$ and $U= (t_i \circ \pi)^{-1}(N)$, therefore $K_i <U$,  $H/K_i \simeq L$ and $U/K_i \simeq N$ is a chief factor of $H$. Let $\beta_i: H/K_i \rightarrow L$ be the isomorphism induced by $t_i \circ \pi$, that is $\beta_i(hK_i):= t_i(\pi(h))$ for all $h \in H$. For all $1 \leq i, j \leq k$ we have $t_i(l) \equiv t_j(l)$ mod $N$ for all $l \in L$. Therefore
\begin{equation} \label{Kj} \beta_i (hK_i) \equiv \beta_j(hK_j) \ \textrm{mod}  \ N , \quad \forall h \in H.\end{equation}

Let $Y := \{ \prescript{g}{} K_i \ | \ g \in G, \ i=1,\ldots, k \}$ and define an equivalence relation $\sim$ on $Y$ as follows: 

We say that $K \sim K'$ if there is an isomorphism $\alpha: H/K \rightarrow H/K'$ such that $\alpha$ induces the identity on $H/U$, i.e. $\alpha(h U/K)= hU/K'$ for each $h \in H$. \medskip

For $1 \leq i, j \leq k$ we can take $\alpha= \beta_j^{-1} \circ \beta_i : H/K_i \rightarrow H/K_j$, which together with (\ref{Kj}) proves that $K_i \sim K_j$.

We claim that the conjugation action of $G$ on $Y$ preserves $\sim$. Indeed if $K \sim K'$ with isomorphism $\alpha$ as above then for any $g \in G$ we have an isomorphism \[ f_{g,K'} \circ \alpha \circ f_{g,K}^{-1}: \ H/ (^g K) \rightarrow H/(^gK'),\] where $f_{g,K}: \ H/ K \rightarrow H/ ^gK$  defined by $f_{g,K}(xK) = \prescript{g}{} x \prescript{g}{}K$ ($ x \in H$) is the isomorphism induced by the conjugation map on $H$. Therefore $^g K \sim \prescript{g}{}K'$ and the claim is proved.

Let $E_1, \ldots, E_l$ be the equivalence classes of $\sim$ on $Y$.  We showed that $K_1 \sim K_2 \sim \cdots \sim K_k$ and without loss of generality we may assume $\{K_1, \ldots, K_k \} \subseteq E_1$. Recall that $Y := \{ \prescript{g}{} K_i \ | \ g \in G, \ i=1,\ldots, k \}$. It follows that $E_1$ intersects nontrivially each orbit of $G$ on $Y$. Therefore the conjugation action of $G$ on $Y$ induces a transitive action of $G$ on the equivalence classes of $\sim$. Let $R_i = \cap_{K \in E_i} K$ and put $R=R_1$. By Proposition \ref{spot} we have $H/R \simeq L_{k'}$ for some $k' \geq k$ and $R \leq M= \cap_{i=1}^k K_i$. 
Further $\cap_{i=1}^l R_i= \cap_{A \in Y} A=\cap_{g \in G} \prescript{g}{} M =T$ as required. 

It remains to prove (\ref{fact}). If $N$ is nonabelian this follows from Proposition \ref{basic} since in that case \[ U/T= U/ (\cap_{K \in Y}K) \simeq  \prod_{K \in Y} U/K\] with $U/R_i \simeq \prod_{K \in E_i} U/K$.

When $N$ is abelian and $K \in Y$ then $H/K \simeq L$ is a split extension isomorphic to $(U/K) \rtimes (H/U)$
and in particular for $K, K' \in Y$ the relation $K \sim K'$ is equivalent to the requirement that the two $H/U$-modules $U/K$ and $U/K'$ are isomorphic. Therefore if $i \not = j$ the $H$-chief factors of $U/R_i$ are not isomorphic (as $H$-modules) to the $H$-chief factors of $U/R_j$. Hence again the Jordan-H\"older theorem implies that $U/(\cap_{i=1}^l R_i)$ being a subdirect product of all $U/R_i$ must in fact be isomorphic to $\prod_{i=1}^l (U/R_i)$. \medskip

Finally the preimage $f^{-1}\left (\{1\}\times  \cdots \times U/R_i \times  \cdots \times \{1\} \right )$ of the direct factor $U/R_i$ of $\prod_{i=1}^l (U/R_i)$ is equal to $\cap_{j \not =i} R_j$ and of course $G$ permutes the $l-1$ element subsets of $\{R_1, \ldots, R_l\}$ in the same way as $G$ permutes $R_1, \ldots R_l$. 

Proposition \ref{spl} is proved. $\square$
 \medskip

Let $R=R_1, R_2, \ldots, R_l$ be the normal subgroups of $H$ from Proposition \ref{spl}.
Let $J_1=N_G(R_1)$ be the normaliser of $R_1$ in $G$. Recall the identity $u(y) \equiv 1$ mod $T$ for all $y \in U$, where  $u(y)=y \cdot \prescript{b_1}{} y^{e_1} \cdots \prescript{b_s}{} y^{e_s}$. Let $b_{i_1}, \ldots, b_{i_t}$ be the subsequence of $b_1, \ldots, b_s$ of elements $b_{i} \in J_1$ and let 
\[ u_0(y): = y \cdot \prescript{b_{i_1}}{} y^{e_{i_1}} \cdots \prescript{b_{i_{t}}}{}y^{e_{i_t}}.\]

We will use the diagonal map $f: U \rightarrow \prod_{i=1}^l (U/R_i)$ from Proposition \ref{spl}, namely $f(y)=(yR_1, \ldots, yR_l)$ for all $y \in U$.

We choose $y_0 \in U$ and let $y \in U$ be an element such that \[ f(y)= (y_0R_1, 1 \ldots 1)  \in \prod_{i=1}^l (U/R_i).\] We have, just as before $f (u(y))= (u_0(y_0), *, \ldots, *)$, and therefore  $u_0(y_0) \in R_1$ for all $y_0 \in U$. 

We have shown that $H/R_1$ is isomorphic to a crown based power $L_{k_1}$ with $k_1 \geq k$ and in particular $d(L_{k_1}) \geq d(L_k)=d(H)$. Since obviously $d(L_{k_1})= d(H/R_1) \leq d(H)$ it follows that $d(H/R_1)=d(H)$.

Further, from (\ref{bd}) we have $d(J_1) \geq d(G) > 6s+s$.

The inequality (\ref{J}) with $J_1$ in place of $J$ gives that Proposition \ref{1/2} (and hence Theorem \ref{iden}) will follow if we prove \[ \frac{d(H)-1}{|J_1:H|} < \alpha(t) \cdot (d(J_1)-1).\]  

Therefore we can replace $G$ with $J_1$, $M$ with $R_1$ and $u$ with $u_0$, and we are reduced to proving the following result.

\begin{prop}  \label{iden1} Let $G$ be a finite group with $d(G)>6s+6$ for some integer $s \geq 1$. Let $H$ be a normal subgroup of $G$, let  $b_i \in G \backslash H$ and $e_i \in \{\pm 1\}$ for $i=1, \ldots, s$. Let $M <H$ be a normal subgroup of $G$ such that $H/M$ is isomorphic to a crown based power $L_k$ with $d(L_k)=d(H)>d(L)$. Let $U/M= \mathrm{soc}(H/M)$ and assume $u(y) \in M$ for all $y \in U$, where   $u(y)= y \cdot \prescript{b_1}{} y^{e_1} \ \cdots \prescript{b_s}{} y^{e_s}$.
Then \[ \frac{d(H)-1}{|G:H|} < \alpha(s) (d(G)-1),\]
where $\alpha(s)=\frac{2s+1}{2s+2}$.
\end{prop}

\subsection{Proof of Proposition \ref{iden1}}

Our argument splits into considering two cases:

\subsubsection{Case 1: the socle $N$ of $L$ is non-abelian.}

First we set up notation for the canonical embedding of $L$ in $Aut(N)$.
For $l \in L$ denote by $\theta(l) \in Aut(N)$ the automorphism  $\theta(l) (x)= \prescript{l}{} x$ ($x \in N$) given by conjugation by $l$. Since $C_L(N)=1$ the map $\theta : L \rightarrow Aut(N)$ is a monomorphism.
Note that $\theta(N)=Inn(N)$, the group of inner automorphisms of $N$. We will denote by $\bar \theta : L/N \rightarrow \theta(L)/\theta(N)$ the induced isomorphism $\bar \theta (lN):=\theta (l) \theta(N)$.
When there is no possibility of confusion we identify $L/N$ with $\theta(L)/\theta(N)$, for example we write $C_{Aut(N)}(L/N)$ below for the centralizer of $\theta(L)/\theta(N)$ in $Aut(N)$. \medskip

We will need a bound for the size of $C_{Aut(N)}(L/N)$. 

Since $N$ is the unique minimal normal subgroup of $L$ it follows that $N \simeq S^m$ a direct product of $m$ copies of some nonabelian simple group $S$. Moreover the conjugation action of $L$ permutes the direct factors of $N$ transitively. Using this it is easy to show that $|C_{Aut(N)}(L/N)| \leq m|N| |Out(S)|$, see proof of Lemma 1 of \cite{VL2}. Now the CFSG gives that $|Out(S)| \leq |S|$ and since $m|S|\leq |S|^m=|N|$ we have 
\begin{equation} \label{estim} |N| \leq |C_{Aut(N)}(L/N)| \leq |N|^2
\end{equation}
\medskip

Let $d=d(G)$ and $d'=d(H)=d(L_k)>d(L)$. Theorem \ref{gen} together with the estimate for $|C_{Aut(N)}(L/N)|$ above gives \begin{equation}  \label{k>} k>P_{L,N}(d'-1) |N|^{d'-1}/|C_{Aut(N)}(L/N)| > \frac{1}{2} |N|^{d'-3}\end{equation} using the fact that $P_{L,N}(d'-1) \geq \frac{53}{90}>1/2$ proved in \cite{DL}.

Fix an isomorphism $\phi: H/M \rightarrow L_k$ and recall that $U/M$ is the socle of $H/M$. Therefore $\phi(U/M)=N^k \leq L_k$, a direct product $N^k$ of $k$ copies of $N$. Let $\pi: L^k \rightarrow L/N$ be the surjection with kernel $N^k$.
Denote by $\rho : H/U \rightarrow L/N$ the isomorphism induced by $\pi \circ \phi$, namely $\rho(hU)=\pi \circ \phi (hM)$ for each $h \in H$.

Let us write $N^k=N_1 \times \cdots \times N_k$, where $N_i \simeq N$ are the minimal normal subgroups of $L_k$ and let $V_i= \phi^{-1}( N_i)$ ($i=1, \ldots, k$). Then $V_1, \ldots, V_k$  are all the minimal normal subgroups of $H/M$ and $U/M= V_1 \times \cdots \times V_k$.

Given a factor $V_i$ of $U/M$ let $P_i/M$ be its complement in $U/M$, i.e. we set $P_i/M:= \prod_{j \not =i} V_j$ for a normal subgroup $P_i < U$ of $H$. Define the isomorphism $\eta_i:H/ P_i \rightarrow L$ as $\eta_i(hP_i)= t_i \circ \phi (hM)$ where $t_i: L_k \rightarrow L$ is the projection in (\ref{ti}). We have $\eta_i(U)=N$. From the definition of $L_k$ it follows that the isomorphism $\bar \eta_i :H/U \rightarrow L/N$ induced by $\eta_i$ is equal to $\rho$ above.

Choose any $b \in G \backslash H$. Let $D$ be the subgroup of $G$ generated by $b$ and $H$.
Since $\{V_1, \ldots, V_k\}$  is the set of the minimal normal subgroups of $H/M$ any automorphism of $H/M$ permutes these $k$ groups among themselves. In particular this applies to the action of $b$ by conjugation on $H/M$.

\begin{prop}\label{fixity} The element $b$ normalizes at most $|N|^{d(D)}$ of the direct factors $V_1, \ldots, V_k$ of $U/M$.
\end{prop}

\textbf{Proof:} Choose $i \in \{1, \ldots, k\}$ and write  $V=V_i$ from now on. Let $P=P_i$ and $\eta= \eta_i$  be as above, i.e. $P/M$ is the complement to $V$ in $U/M$ and $\eta: H/P \rightarrow L$ is the isomorphism associated to the projection $t_i \circ \phi : H \rightarrow L$.

Suppose that  $\prescript{b}{} V=V$. This is equivalent to $^b (P/M)=P/M$ i.e. $P$ is normalized by $D=\langle b, H \rangle$.

The conjugation action of
$D$ on $U/P$ together with $\eta$ induce a homomorphism $f: D \rightarrow Aut(N)$ defined by \[ f(g) (x)= \eta ( \prescript{g}{} \eta^{-1} (x)) \quad  \forall g \in D, \forall x \in N.\]  

Note that if $h \in H$ then $f(h)$ is the conjugation by $\eta(hP)$ on $N$ and therefore $f(h)=\theta(\eta(hP))$. In particular $f(H)=\theta(L) \leq Aut(N)$ and $f(U)=\theta(N)=Inn(N)$.

The homomorphism $\bar f : H/U \rightarrow \theta(L)/\theta(N)$ induced by $f$
is given by 
\[\bar f(hU):= f(h)\theta (N)=\theta(\eta(hP)N)\theta (N)= \bar \theta (\eta (hP)N)=\bar \theta \circ \bar \eta(hU) \] for all $h \in H$. Therefore $\bar f=\bar \theta \circ \bar \eta=\bar \theta \circ \rho$ and does not depend on the choice of $V=V_i$. Observe that $f$ determines $P$ (and hence $V$) since $P= \ker f \cap H$.

Define the set $\mathfrak F$ to consist of all homomorphisms $f: D \rightarrow Aut(N)$ such that $f(H)=\theta(L)$, $f(U)=\theta(N)=Inn(N)$ and $\bar f = \bar \theta \circ \rho$ where $\bar f: H/U \rightarrow \theta (L)/\theta(N)$ is the isomorphism induced by $f$.

We summarize the above discussion in the following.

\begin{prop} \label{F} The number of direct factors $V_i$ of $U/M$ normalized by $b$ is bounded above by the number of different subgroups $H \cap \ker f$ where $f \in \mathfrak F$. \end{prop}

We will prove that $|\mathfrak F | \leq |C_{Aut(N)} (L/N)| \cdot |N|^{d(D)}$.

Let $f \in \mathfrak F$ and let $\sigma:= f(b) \in Aut(N)$. We first estimate the possibilities for $\sigma$.
For $h \in H$ we have \[ \bar \theta \circ \rho (\prescript{b}{} hU)= f(\prescript{b}{} h) \theta(N) = \prescript{\sigma}{} f(h)\theta(N)=\prescript{\sigma}{} (\bar \theta \rho (hU)).\]

Thus the action of $\sigma$ by conjugation on $\theta(L)/\theta(N)$ is uniquely determined by $\rho$ and the conjugation action of $b$ on $H/U$. So any two choices $\sigma$ and $\sigma'$ for $f(b)$ satisfy that $\sigma^{-1} \sigma'$ centralizes $\theta(L)/\theta(N)$ . Hence there are at most $|C_{Aut(N)}(L/N)|\leq  |N|^2$ possibilities for $\sigma=f(b)$ in $Aut(N)$.

Let $Q:=f(D) \leq Aut(N)$ be the image of $f$ in $Aut(N)$. Observe that $Q=\langle \theta(L), \sigma \rangle$ since $D=\langle H, b\rangle$ and $f(H)=\theta(L)$. Therefore $\sigma$ uniquely determines $Q$.

The induced map $\tilde f : D/U \rightarrow Q/\theta(N) \leq Out(N)$ defined by $\tilde f (gU)= f(g) \theta(N)$ for all $g \in D$, is also uniquely determined by $\sigma=f(b)$. Indeed, $D/U$ is generated by $H/U$ together with $bU \in G/U$ as a subgroup of $G/U$. The restriction $\tilde f|_{H/U}$ of $\tilde f$ to $H/U$ equals $\bar f=\bar \theta \circ \rho$, while of course $\tilde f(bU)= \sigma \theta (N)$ is determined by $\sigma$.

Therefore any choice of $\sigma=f(b)$ uniquely determines the group $f(D)=Q \leq Aut(N)$ and $\tilde f: D/U \rightarrow Q/\theta(N)$. 

Next we estimate the possibilities for $f$ given $Q$ and $\tilde f$. Let $r=d(D)$ and choose a generating set $g_1, \ldots, g_r$ of $D$. We count how many possibilities are there for the images $f(g_1), \ldots, f(g_r)$ in $Q$. Since $f(g_i)\theta(N)= \tilde f (g_iU)$ is a uniquely determined coset of $Q/\theta(N)$  we obtain that there are at most $|N|$ choices for each $f(g_i) \in Q$. Thus for any chosen $\sigma =f(b)$ there at most $|N|^{r}$ possibilities for $f$, and hence $|\mathfrak F | \leq |C_{Aut(N)} (L/N)||N|^{r}$. 

Next we define an action of $C_{Aut(N)}(L/N)$ on $\mathfrak F$ by

\[ \prescript{\gamma}{} f (g):= \gamma \circ f(g) \circ \gamma^{-1} \quad \textrm{for all } f \in \mathfrak F, \ g \in D, \ \gamma \in C_{Aut(N)}(L/N). \]

We claim that $C_{Aut(N)}(L/N)$ acts semi-regularly on $\mathfrak F$. Indeed, suppose $f= \prescript{\gamma}{} f$ i.e. $f(g)= \gamma \circ f(g) \circ \gamma^{-1}$ for all $g \in D$. In particular, by letting $g$ range over $U$ and using $f(U)=\theta(N)$ we deduce that $\theta(x) \circ \gamma = \gamma \circ \theta(x)$ for all $x \in N$. This means that $ \prescript{x}{} \gamma (y)= \gamma( \prescript{x}{} y)=\prescript{\gamma(x)}{}\gamma(y)$ for all $x,y \in N$.  Therefore  $x^{-1}\gamma(x) \in N$ centralizes all elements of $N$ and so $x= \gamma(x)$ since $N$ is a nonabelian minimal normal subgroup of $L$. Thus $\gamma=1$ and hence the action of $C_{Aut(N)}(L/N)$ on $\mathfrak F$ is indeed semi-regular.

Finally observe that $\ker f= \ker (\prescript{\gamma}{} f)$. Therefore by Proposition \ref{F} the number of factors $V_i$ normalized by $b$ is bounded above by the number of orbits of $C_{Aut(N)}(L/N)$ on $\mathfrak F$, which is at most $|N|^r$.
Proposition \ref{fixity} follows.  
$\square$
\medskip

We will apply Proposition \ref{fixity} with $b=b_i \in G\backslash H$ from the statement of Proposition \ref{iden1}. Note that $D=\langle H,b \rangle \not = H$.

Recall that $d'=d(H)$ and $d=d(G)$. 

We have $r=d(D) \leq  (d-1)|G:D|+1$ and from $D \not = H$ we have $|G:D| \leq |G:H|/2$, whence $r \leq \frac{(d-1)|G:H|}{2}+1$.

Proposition \ref{fixity} gives that each $b_i$ can normalize at most $A:=|N|^{|G:H|\frac{d-1}{2}+1}$ of the factors $V_i$ of $U/M$.

Assume, for the sake of contradiction, that $d'=d(H)\geq  \frac{3}{4} |G:H|(d-1)+1$ and recall that $k > \frac{1}{2} |N|^{d'-3}$ from (\ref{k>}).  Therefore 
\[ k >  \frac{1}{2}|N|^{\frac{3}{4}|G:H|(d -1)-2}= \frac{1}{2} A |N|^{\frac{1}{4}|G:H|(d-1)-3} . \]
Since $|G:H| \geq 2$ and $|N| \geq 60$ we have \[ k >  \frac{1}{2} A \ |N|^{\frac{d-1}{2}-3} \geq  \frac{1}{2}A \ 60^{\frac{d-7}{2}}.\]
Using that $60^a>a \ln 60>4a$ for $a>0$ and $d=d(G) >6s+6>s+7$ we obtain $\frac{1}{2} 60^{\frac{d-7}{2}}>d-7 > s$ and hence $k>As$. 

It follows that there is a direct factor, say $V_1$ of $U/M =V_1 \times \cdots \times V_k$ such that $\ ^{b_i} V_1 \not = V_1$ for each $i=1, \ldots, s$. 

Recall the identity $u(y)= y \cdot \prescript{b_1}{} y^{e_1} \ \cdots \prescript{b_s}{} y^{e_s} \equiv 1$ mod $M$ for all $y \in U$. Choose $y \in V_1\backslash\{1\}$. It follows that the projection of $u(y)$ onto $V_1$ is equal to $y \not =1$, hence $u(y) \not = 1$ mod $M$, contradiction. Therefore \[ d(H) < \frac{3}{4}|G:H|(d-1)+1 \] and Proposition \ref{iden1} follows since $\alpha(s) \geq \frac{3}{4}$.   

\subsubsection{Case 2: the socle $N$ of $L$ is abelian.}

Let $p$ be the exponent of the abelian group $N$, then $p$ is a prime and $N$ is a simple $\mathbb F_p(L/N)$-module.

Recall that $U/M$ is the socle of $H/M$. Thus $U$ is a normal subgroup of $G$ and we define $\Gamma=G/U$ and $\Delta= H/U$. We will consider $U/M$ as a left module for the action of $\Gamma$ and $\Delta$ by conjugation.

Fix an isomorphism $\phi: H/M \rightarrow L_k$ with $\phi(U/M)=N^k=soc(L_k)$. Then $\phi$ induces an isomorphism  $\bar \phi: \Delta \rightarrow L_k/N^k \simeq L/N$ given by $\bar \phi (hU)=\phi(hM) N^k$ for all $h \in H$. From now on we will consider $N$ as $\mathbb F_p \Delta$-module via $\bar \phi$. 

The isomorphism between $U/M$ and  $\phi(U/M)= soc(L_k)=N^k$ gives that
$U/M$ is a semisimple $\mathbb F_p \Delta$-module which is a direct sum of $k$ isomorphic copies of $N$. 

Let $B \leq U/M$ be any simple $\mathbb F_p \Delta$-submodule of $U/M$. Then $B$ is isomorphic to $N$. Let $g \in \Gamma$. The submodule $g^{-1}B \leq U/M$ is also a simple $\mathbb F_p \Delta$-submodule of $U/M$ and hence is isomorphic to $B$ and $N$. Let $\sigma: B \rightarrow B$ be the linear map such that $x \mapsto g^{-1} \sigma(x)$ is an isomorphism of $\mathbb F_p \Delta$-modules between $B$ and $g^{-1}B$. This means that $g^{-1} \sigma(h x)=h g^{-1}  \sigma(x)$ for all $x \in B$ and $h \in \Delta$, which is equivalent to $\sigma( h x)= \prescript{g}{} h \cdot \sigma(x)$. Using that $B$ and $N$ are isomorphic $\mathbb F_p \Delta$-modules we deduce the following.

\begin{prop} \label{twist} Let $g \in \Gamma$. There is a bijective linear map $\sigma=\sigma_g: N \rightarrow N$ such that $\sigma (hx)= \prescript{g}{}h  \ \sigma (x)$ for all $x \in N$ and $h \in \Delta$.
\end{prop}  

For a $\mathbb F_p \Delta$-module $A$ define  \[ \mathfrak R_N(A)= \cap \{\ker f \ | \ f \in Hom_\Delta (A, N) \},\] where $f$ ranges over all $\mathbb F_p \Delta$-homomorphisms from $A$ to $N$.
 
Thus $A/ \mathfrak R _N(A)$ is the largest semisimple quotient of $A$ whose simple factors are isomorphic to $N$. Suppose that $A$ is in addition a $\mathbb F_p\Gamma$-module with a $\mathbb F_p \Delta$-submodule $B$. Let $f \in Hom_\Delta(B,N)$ and $g \in \Gamma$.  Define $^g f: gB \rightarrow N$ by $^gf (x):=\sigma_g \circ f(g^{-1}x)$ for all $x \in gB$. Note that $^gf \in Hom_\Delta (gB,N)$ and $\ker \prescript{g}{}f = g \ker f$. Therefore $\mathfrak R_N(gB) \leq g \mathfrak R_N(B)$ and by reversing the role of $B$ and $gB$ we get 
\begin{equation}\label{rn}\mathfrak R_N(gB) = g \mathfrak R_N(B). \end{equation}
 In particular $\mathfrak R_N(A)$ is a $\mathbb F_p \Gamma$-submodule of $A$.

We would need to know more about the action of $\Gamma$ on $U/M$.
 
Let $\mathcal G$ be the free group on $d=d(G)$ generators $x_1, \ldots, x_d$ and fix a surjective homomorphism $\pi: F \rightarrow G/M$. Let $\mathcal H, \mathcal U$ be the preimages of $H/M$ and $U/M$ respectively and let  $\mathcal M = \ker \pi$. We will identify $\Gamma= G/U$ and $\Delta=H/U$ with $\mathcal G/\mathcal U$ and $\mathcal H/\mathcal U$ via $\pi$. Thus the semisimple $\mathbb F_p \Gamma$-module $U/M$ is isomorphic to the 
$\mathbb F_p(\mathcal G/\mathcal U)$-module $\mathcal U/ \mathcal M$.  In turn $\mathcal U/\mathcal M$ is a quotient of the $\mathbb F_p \Gamma$-module $\mathcal U_p^{ab}=\mathcal U/\mathcal U^p [\mathcal U, \mathcal U]$, the mod-$p$ relation module of the presentation \[ 1 \rightarrow \mathcal U \rightarrow \mathcal G \rightarrow \Gamma \rightarrow 1. \]  

We will denote by $\Phi(\mathcal U)=\mathcal U^p[\mathcal U, \mathcal U]$ the $p$-Frattini subgroup of $\mathcal U$.

The cellular chain complex of the Cayley graph of $\Gamma$ with respect to the generators $\pi(x_1), \ldots, \pi(x_d)$ gives rise to the exact sequence of $\mathbb F_p \Gamma$-modules
\begin{equation} 0 \rightarrow \mathcal U_p^{ab} \stackrel{\rho}{\rightarrow} (\mathbb F_p \Gamma)^d \rightarrow \mathbb F_p \Gamma \rightarrow \mathbb F_p \rightarrow 0.
\end{equation}

The injection $\rho: \mathcal U_p^{ab} \rightarrow (\mathbb F_p \Gamma)^d $ can be described explicitly by
\begin{equation} \label{rho} \rho(u \Phi(\mathcal U))= \left ( \overline{\frac{\partial u}{\partial x_1}}, \ldots, \overline{\frac{\partial u}{\partial x_d}} \right), \quad  \forall u \in \mathcal U.\end{equation}
 Here $\frac{\partial u}{\partial x_i} \in \mathbb F_p \mathcal G$ is the Fox derivative of $u$ with respect to $x_i$, and $\bar y$ denotes the image of $y \in \mathbb F_p \mathcal G$ under the reduction $\mathbb F_p \mathcal G \rightarrow \mathbb F_p \Gamma$. See \cite{B}, \S I.5, Proposition (5.4) and Exercise 3(d) for proofs.
 
Let $\mathcal M_0= \mathcal M/ \Phi(\mathcal U)$ be the image of $\mathcal M$ in $\mathcal U_p^{ab}$, thus $U/M$, $\mathcal U/\mathcal M $ and $\mathcal U_p^{ab} / \mathcal M_0$ are all isomorphic as $\mathbb F_p \Gamma$-modules.

The following result can be deduced from the fact that for any $\mathbb F_p \Delta$-module $S$ we have \[ H^2(\Delta,S)= \mathrm{Coker} \left ( Hom_\Delta((\mathbb F_p \Gamma)^d,S) \rightarrow Hom_\Delta(\mathcal U^{ab},S)  \right ) \] and \cite{B}, \S IV.2, Exercise 4. For completeness we give a proof of it in section \ref{prooflift}.

\begin{prop} \label{lift} Let 
$\mathcal K$ be a normal subgroup of $\mathcal H$ with $\Phi(\mathcal U) \leq \mathcal K \leq \mathcal U$ and let $S:= \mathcal U/\mathcal K$ considered as $\mathbb F_p \Delta$-submodule under conjugation by $\mathcal H$. Let $f: \mathcal U_p^{ab} \rightarrow S$ be the associated epimorphism of $\mathbb F_p \Delta$-modules, namely $f(u \Phi(\mathcal U))=u\mathcal K$ for all $u \in \mathcal U$. Then $\mathcal U/\mathcal K$ has a complement in $\mathcal H/\mathcal K$ if and only if there is a $\mathbb F_p \Delta$-module homomorphism $\theta: (\mathbb F_p \Gamma)^d \rightarrow S$ such that $f= \theta \circ \rho$.
\end{prop}

We claim that $\mathcal M$  is an intersection of groups $\mathcal K_i  \leq \mathcal U$ for $i=1, \ldots, k$ such that $\mathcal U/\mathcal K_i$ is $\mathbb F_p \Delta$-isomorphic to $N$ and has a complement in $\mathcal H/\mathcal K_i$. In view of the isomorphisms $\mathcal H/\mathcal M \rightarrow H/M \rightarrow L_k$
it is sufficent to prove the corresponding statement for $\{0\}$ in $L_k$ and take preimages. For $i=1, \ldots, k$ let $t_i: L_k \rightarrow L$ be the projection defined in (\ref{ti}). Then $\cap_{i=1}^k \ker t_i=\{0\}$ and $L_k/\ker t_i \simeq L$ which splits over $N$.
This proves the claim.

The groups $\mathcal K_i$ satisfy the requirements of Proposition \ref{lift}, hence if $f_i : \mathcal U_p^{ab} \rightarrow N$ are epimorphisms with  $\ker f_i = \mathcal K_i/\Phi(\mathcal U)$ there exist $\mathbb F_p \Delta$-epimorphisms $\theta_i: (\mathbb F_p \Gamma)^d \rightarrow N$ such that $f_i= \theta_i \circ \rho$.

Let $Q:= \mathfrak R_N((\mathbb F_p\Gamma)^d)$ and note that $Q$ is a $\mathbb F_p \Gamma$ submodule of $(\mathbb F_p\Gamma)^d$. 

Now we prove $Q \cap \rho(\mathcal U_p^{ab}) \leq \rho(\mathcal M_0)$. Let $x \in \mathcal U_p^{ab}$ be such that $\rho(x) \in Q$. In particular $f_i(x)=\theta_i (\rho (x))=0$ for $i=1, \ldots, k$ and hence  $x$ belongs to \[ \bigcap_{i=1}^k \ker f_i =\bigcap_{i=1}^k \mathcal K_i/\Phi(\mathcal U)= \mathcal M/ \Phi(\mathcal U)=\mathcal M_0. \]

Therefore $Q \cap \rho(\mathcal U_p^{ab}) \leq \rho(\mathcal M_0)$ which gives $\rho(\mathcal U_p^{ab}) \cap (\rho(\mathcal M_0) +Q)= \rho(\mathcal M_0)$. In particular we have the following isomorphisms of $\mathbb F_p \Gamma$ modules 
\[ \frac{U}{M} \simeq \frac{\mathcal U}{\mathcal M} \simeq \frac{\mathcal U_p^{ab}}{\mathcal M_0} \simeq \frac{\rho(\mathcal U_p^{ab})}{\rho(\mathcal M_0)} \simeq  \frac{\rho(\mathcal U_p^{ab})+Q}{\rho(\mathcal M_0)+Q}. \] 

Let $\bar b_i=b_iU$ be the images of $b_i \in G$ in $\Gamma=G/U$ and let $W:=\frac{(\mathbb F_p \Gamma)^d}{Q}$. The condition $u(y) \in M$ for all $y \in U$  in Proposition \ref{iden1} is equivalent to
\[ (1+e_1\bar b_1 + \cdots +e_s\bar b_s) \frac{U}{M}= \{0\}. \] Therefore $1+e_1\bar b_1 + \cdots +e_s\bar b_s$ annihilates the sub-quotient $\frac{\rho(\mathcal U_p^{ab})+Q}{\rho(\mathcal M_0)+Q}$ of $W$.  

Recall a basic result from linear algebra 
\begin{prop} Let $V$ be a vector space over a field $\mathbb F$ and let $T: V \rightarrow V$ be a linear transformation with a pair $X \geq Y$ of $T$-invariant subspaces such that $T (X/Y)=\{0\}$ i.e. $T(X) \subseteq Y$. Then $\dim_{\mathbb F} \ker T \geq \dim_{\mathbb F} X/Y$. 
\end{prop}
In particular we deduce 

\begin{equation}\label{dim} \dim_{\mathbb F_p} \ker (1+ e_1\bar b_1 + \cdots + e_s\bar b_s)|_W  \geq \dim_{\mathbb F_p} U/M= k \dim_{\mathbb F_p} N.\end{equation}

We now investigate the structure of $W$. Let $m:=|\Gamma: \Delta|$ and choose coset representatives $\gamma_1, \gamma_2, \ldots, \gamma_m$ for $\Gamma/\Delta$ in $\Gamma$.
We can write $\mathbb F_p \Gamma= \oplus_{i=1}^m \gamma_i (\mathbb F_p \Delta)$ as a direct sum where each factor $\gamma_i (\mathbb F_p \Delta)= (\mathbb F_p \Delta) \gamma_i$ is a free $\mathbb F_p \Delta$-module.
In particular  (\ref{rn}) gives \[ \mathfrak R_N(\mathbb F_p \Gamma)= \bigoplus_{i=1}^m \mathfrak R_N (\gamma_i (\mathbb F_p \Delta))= \bigoplus_{i=1}^m \gamma_i \mathfrak R_N (\mathbb F_p \Delta).\]
Therefore
\[ W=\frac{(\mathbb F_p \Gamma)^d}{\mathfrak R_N((\mathbb F_p \Gamma)^d)}= \left( \bigoplus_{i=1}^m \frac{\gamma_i \mathbb F_p \Delta}{\gamma_i \mathfrak R_N(\mathbb F_p \Delta)}\right)^d. \]

Let $V_i=\left (\frac{\gamma_i \mathbb F_p \Delta}{\gamma_i \mathfrak R_N(\mathbb F_p \Delta)}\right)^d$ so that $W= V_1 \oplus \cdots \oplus V_m$. If $\gamma \in \Gamma$ then the action of $\gamma$ on $W$ permutes the spaces $V_1, \ldots, V_m$ in the same way as $\gamma$ acts on $\Gamma/\Delta$, namely $\gamma V_i=V_j$ where $\gamma \gamma_i \Delta= \gamma_j \Delta$.

We need the following elementary result.
\begin{prop} \label{linalg} Let $V=V_1 \oplus \cdots \oplus V_m$ be a vector space over a field $\mathbb F$ which decomposes as a  direct sum of its subspaces $V_1, \ldots, V_m$. Let $s \in \mathbb N$  and for $i=1, \ldots, s$ let $T_i: V \rightarrow V$  be bijective linear transformations preserving the above direct sum decomposition. Assume that for all $i=1, \ldots, s$ and $j=1, \ldots, m$ we have $T_i(V_j) \in \{V_1, \ldots, V_m\} \backslash \{V_j\}$. Then \[ \dim_{\mathbb F} \ker ( 1+T_1 + \cdots + T_s) \leq \frac{s}{s+1} \dim_{\mathbb F} V.\]
\end{prop}
\textbf{Proof:} Let $\mathfrak G \leq GL(V)$ be the subgroup of $GL(V)$ generated by $T_1, \ldots, T_s$. Let $E_1, \ldots E_r$ be the orbits of $\mathfrak G$ acting on the direct summands $\{V_1, \ldots, V_m\}$ and for $j=1,\ldots, r$ define $\mathfrak V_j= \oplus_{V_i \in E_j} V_i$. Since each $\mathfrak V_j$ is invariant under $\mathfrak G$ and satisfies the same hypothesis as $V$, it will be sufficient to prove the claimed inequality for each $\mathfrak V_j$ in place of $V$, and then add them for $j=1, \ldots, r$. Therefore without loss of generality we  may assume that  $\mathfrak G$ acts transitively on $V_1, \ldots, V_m$. In particular, since $\dim V_i =\dim g(V_i)$ for each $V_i$ and each $g \in \mathfrak G$, we have $\dim V_1 = \cdots = \dim V_m$.

Let $t = \lceil \frac{m}{s+1} \rceil$. We construct inductively a sequence $V_{a_1}, \ldots, V_{a_t}$ as follows: Put $V_{a_1}=V_1$ and for $i=2, \ldots ,t$ choose $V_{a_i}$ to be any element in 
\[  \{V_1, \ldots, V_m\} \backslash \left (\bigcup_{j=1}^{i-1} \left \{ V_{a_j}, T_1(V_{a_j}), \ldots , T_s(V_{a_j}) \right \} \right) .\]   
Such $V_{a_i}$ exists as long as $(i-1)(s+1) < m$ i.e. $i \leq t$. 

Let $\pi_i: V \rightarrow V_i$ be the projection onto $V_i$.  The choice of $\{V_{a_i}\}_{i=1}^t$ ensures that for all pairs $1 \leq j <i \leq t$ \[ \pi_{a_i} \circ ( 1+T_1 + \cdots +T_s)(V_{a_j}) =\{0\}, \] 
while the condition $T_r(V_{a_i}) \not = V_{a_i}$ for $r=1, \ldots, s$ ensures that \[ \pi_{a_i} \circ (1+T_1 + \cdots T_s)|_{V_{a_i}}= \mathrm{1}_{V_{a_i}}.\]  

It follows that the $t$ subspaces $\{(1+T_1 + \cdots +T_s)(V_{a_i})\}_{i=1}^t$ generate their direct sum in $V$. Moreover $\dim (1+T_1 + \cdots +T_s)(V_{a_i})= \dim V_{a_i}= \dim V_1$ for all $i=1, \ldots, t$.

Therefore \[\dim (1+T_1 + \cdots T_s)(V) \geq t \dim V_1 \geq \frac{m}{s+1} \dim V_1= \frac{1}{s+1} \dim V\] and Proposition \ref{linalg} follows from the rank-nullity theorem. $\square$
\medskip

Let $J (\mathbb F_p \Delta)$ be the Jacobson ideal of  $\mathbb F_p \Delta$.
By the theory of semisimple algebras the multiplicity of $N$ in $\mathbb F_p \Delta/ J (\mathbb F_p \Delta)$ is equal to $c:=\dim_E N$, where $E= End_{\Delta}(N)$ is a finite extension field of $\mathbb F_p$. Therefore 
$\frac{\mathbb F_p \Delta}{\mathfrak R_N(\mathbb F_p \Delta)} \simeq N^c$  and hence \[ \dim_{\mathbb F_p} W= d\cdot m \cdot c \cdot \dim_{\mathbb F_p} N, \] where $d=d(G)$ and $m=|G:H|=|\Gamma : \Delta|$.

Apply Proposition \ref{linalg} with $V= W$, $V_i=  \left ( \frac{\gamma_i \mathbb F_p \Delta}{\gamma_i \mathfrak R_N(\mathbb F_p \Delta)}\right)^d$ and $T_j(x)=e_j\bar b_j x$ for $x \in W$  and $j=1, \ldots, s$. Since $\bar b_j \not \in \Delta$  we have $ \bar b_j \gamma_i \Delta \not = \gamma_i \Delta$ for any coset $\gamma_i\Delta$ and hence $\bar b_j V_i \not = V_i$. We deduce 
\[ \dim_{\mathbb F_p} \ker (1 + e_1\bar b_1 + \cdots + e_s\bar b_s)|_W \leq \frac{s}{s+1} \dim_{\mathbb F_p} W=\frac{s}{s+1} d \cdot m \cdot c \dim_{\mathbb F_p} N.\] From (\ref{dim}) we have $\dim_{\mathbb F_p} \ker (1 + \bar b_1 + \cdots + \bar b_s)|_W \geq k \dim_{\mathbb F_p} N$ and therefore  $k \leq \frac{sdmc}{s+1}$ On the other hand, if we set $d'=d(H)=d(L_k)$, Theorem \ref{gen} gives \[ k > (d'-2)c - \dim_E H^1(L/N,N).\] Moreover, since $N$ is a faithful $L/N$-module we have $|H^1(L/N,N)| < |N|$ by \cite{AG} Theorem A, and therefore $k >(d'-2)c-c=(d'-3)c$. Hence \[ (d'-3)c < k \leq \frac{s  d m c}{s+1} \] and so $d'<\frac{sdm}{s+1}+3$. Now 
\[\frac{sdm}{s+1}+3= \frac{(2s+1)(d-1)m}{2s+2} +1 + C ,\] where \[ C=\frac{2s+1}{2s+2}m+2 - \frac{dm}{2s+2}< m+2-3m \leq 0 \]
since $d=d(G)>6(s+1)$. 
Therefore \[ d'=d(H) <  \frac{sdm}{s+1}+3 < \frac{(2s+1)(d-1)m}{2s+2} +1 \] and recalling that $m= |\Gamma: \Delta|=|G:H|$ we obtain
\[\frac{d(H)-1}{|G:H|}= \frac{d'-1}{m} < \frac{2s+1}{2s+2} (d-1) = \alpha(s)  (d(G)-1). \]
The proof of Proposition \ref{iden1} is complete. 

\subsection{Proof of Proposition \ref{lift}.} \label{prooflift}

First we show that it suffices to prove Proposition \ref{lift} in the special case when $\mathcal H = \mathcal G$.
Let $d=d(\mathcal G)$ and let $B_d$ be the graph with a single vertex and $d$ edges, namely the bouquet of $d$ circles. We identify $\mathcal G$ with $\pi_1(B_d)$, the fundamental group of $B_d$.
Let $\tilde{B_d}$ be the universal cover of $B_d$, this is just the Cayley graph of $\mathcal G$. The cellular chain complex $R_1:= C_*(\tilde{B_d}, \mathbb F_p)$ of $\tilde{B_d}$ together with the action of $\mathcal G$ by translation on $\tilde B_d$ results in the free $\mathbb F_p \mathcal G$- resolution of $\mathbb F_p$.
\[ R_1: \quad 0 \rightarrow (\mathbb F_p \mathcal G)^d \rightarrow \mathbb F_p \mathcal G \rightarrow \mathbb F_p \rightarrow 0. \] 
Let $Z$ be the quotient graph of $\tilde B_d$ by the action of $\mathcal H \leq \mathcal G$. Then $Z$ is a cover of $B_d$ of degree $|\mathcal G: \mathcal H|$ with fundamental group $\pi_1(Z)=\mathcal H$. In fact $Z$ is the Cayley graph of $\mathcal G/\mathcal H$ with respect to the image of the generating set of $\mathcal G$. 

By contracting a maximal tree in $Z$ to a point we obtain a homotopy equivalence \[ \alpha: Z \rightarrow  B_{d'}\] of $Z$ with the bouquet of $d'$ circles $B_{d'}$, where $d'=d(\mathcal H)=(d-1)|\mathcal G : \mathcal H| +1$. The homotopy equivalence $\alpha$ induces an isomorphism of fundamental groups $\alpha_\pi: \mathcal H =\pi_1(Z) \rightarrow \pi_1(B_{d'})$ and we will identify $\mathcal H$ with $\pi_1(B_{d'})$ via $\alpha_\pi$ from now on. Moreover, $\alpha$ induces a homotopy equivalence of chain complexes of $\mathbb F_p \mathcal H$-modules between $R_1=C_*(\tilde Z,\mathbb F_p)$ and the cellular chain complex $R_2:=C_*(\tilde{B}_{d'}, \mathbb F_p)$ of the universal cover $\tilde B _{d'}$ of $B_{d'}$, namely
\[ R_2: \quad 0 \rightarrow (\mathbb F_p \mathcal H)^{d'} \rightarrow \mathbb F_p \mathcal H \rightarrow \mathbb F_p \rightarrow 0. \]   

 Specifically, there are chain maps of complexes of $\mathbb F_p \mathcal H$-modules $\tilde \alpha_* : R_1 \rightarrow R_2$ and $\tilde \beta_*: R_2 \rightarrow R_1$ such that the composition $\tilde \alpha_* \circ \tilde \beta_*$ (resp. $\tilde \beta_* \circ \tilde \alpha_*$) is homotopic to the identity map on $R2$ (resp. $R1$). Passing to $\mathcal U$-coinvariants by applying $\mathbb F_p \otimes_{\mathbb F_p\mathcal U }  - $ (i.e. factoring out the action of $\mathcal U$), we deduce that $\tilde \alpha_*, \tilde \beta_*$ induce homotopy equivalence between the complexes of $\mathbb F_p \Delta$-modules 
 \[ (\mathbb F_p \Gamma)^d  \stackrel{\delta'}{\longrightarrow} \mathbb F_p \Gamma \rightarrow \mathbb F_p \rightarrow 0\] and
\[ (\mathbb F_p \Delta)^{d'} \stackrel{\delta''}{\longrightarrow} \mathbb F_p \Delta \rightarrow \mathbb F_p \rightarrow 0 \]
In particular there are $\mathbb F_p \Delta$-homomorphisms in degree 1:
$\bar \alpha : (\mathbb F_p \Gamma)^d \rightarrow (\mathbb F_p \Delta)^{d'}$ and $\bar \beta : (\mathbb F_p \Delta)^{d'} \rightarrow (\mathbb F_p \Gamma)^d$. which induce mutually inverse bijections between $X_1= \ker \delta'$ and $X_2 =\ker \delta''$. The canonical isomorphisms $\rho_i: \mathcal U_p^{ab} \rightarrow  X_i$ from (\ref{rho}) give $\bar \alpha|_{X_1}=\rho_2 \circ \rho_1^{-1}$ and $\bar \beta|_{X_2}= \rho_1 \circ \rho_2^{-1}$ since the homotopy $\alpha$ induces the identity on the fundamental group $\mathcal H$ with our identification. 

Let $S$ be a $\mathbb F_p \Delta$-module and let $f: \mathcal U_p^{ab} \rightarrow S$ be a $\mathbb F_p \Delta$-homomorphism. If $f=\theta \circ \rho_1$  for some $\theta \in  Hom_{\Delta} ((\mathbb F_p \Gamma)^d ,S)$ then $f= \zeta \circ \rho_2$, where $\zeta:= \theta \circ \bar \beta$ belongs to  $Hom_{\Delta}((\mathbb F_p \Delta)^{d'}, S)$. Conversely, if $f= \zeta \circ \rho_2$ for  $\zeta \in Hom_{\Delta}((\mathbb F_p \Delta)^{d'}, S)$ then $f=\theta \circ \rho_1$ with $\theta:= \zeta \circ \bar \alpha$ in $Hom_{\Delta} ((\mathbb F_p \Gamma)^d ,S)$.

Hence the validity of Proposition \ref{lift} depends only on the presentation \[ 1 \rightarrow \mathcal U \rightarrow \mathcal H \rightarrow \Delta \rightarrow 1\] without reference to $\mathcal G$, and we may assume from now on that $\mathcal G =\mathcal H$ and $\Gamma=\Delta$.

Let us restate the notation of Proposition \ref{lift} in this setting.

Let $d=d(\mathcal G)$ and fix a generating set $x_1, \ldots, x_{d}$ of $\mathcal G$. The isomorphism $\rho: \mathcal U_p^{ab} \rightarrow (\mathbb F_p \Gamma)^{d} $ given by (\ref{rho}) completes the exact sequence  of $\mathbb F_p \Gamma$-modules \[ 0\rightarrow \mathcal U_p^{ab} \stackrel{\rho}{\longrightarrow}  (\mathbb F_p \Gamma)^{d} \rightarrow  \mathbb F_p \Gamma \rightarrow \mathbb F_p \rightarrow 0 \] associated to the presentation $1 \rightarrow \mathcal U \rightarrow \mathcal G \rightarrow \Gamma \rightarrow 1$.

Let $\mathcal K$ be a normal subgroup of $\mathcal G$ with $\Phi(\mathcal U) \leq \mathcal K \leq \mathcal U$. Let $S$ be the $\mathbb F_p \Gamma$-module $\mathcal U/\mathcal K$ and  let $f \in Hom_\Gamma(\mathcal U_p^{ab}, S)$ be the map $u \Phi(\mathcal U) \mapsto u\mathcal K$  ($u \in \mathcal U$). We need to show that $f$ lifts to a $\mathbb F_p \Gamma$-homomorphism $\theta: (\mathbb F_p \Gamma)^{d} \rightarrow S$ such that $f= \theta \circ \rho$ if and only if $\mathcal U/\mathcal K$ is complemented in $\mathcal G/\mathcal K$.

We have the following basic result.
\begin{prop} \label{split}
The group $\mathcal U/ \mathcal K$ is complemented in $\mathcal G /\mathcal K$ if and only if there is a homomorphism $\phi: \mathcal G \rightarrow \frac{\mathcal G}{\mathcal K}$ such that 

P1. $\phi(y) \equiv y\mathcal K$  mod $\mathcal U/\mathcal K$ for each $y \in \mathcal G$, and 

P2. $\mathcal U \subseteq \ker \phi$. 
\end{prop}
\textbf{Proof:}
Suppose there is a homomorphism $\phi: \mathcal G \rightarrow \mathcal G/\mathcal K$ with properties P1 and P2.  Define $\Theta:=\phi(\mathcal G) \leq \mathcal G/\mathcal K$.  Suppose $\phi(x) \in (\mathcal U/\mathcal K) \cap \Theta$ for some $x \in \mathcal G$. From P1 $ x\mathcal K \equiv \phi(x)$ mod  $\mathcal U/\mathcal K$ and therefore $x \in \mathcal U$. By P2 we have $\phi(x)=1$ and we have proved $\Theta \cap (\mathcal U/\mathcal K)=1_{\mathcal G/\mathcal K}$. In addition P1 gives $\Theta \cdot (\mathcal U/\mathcal K)= \mathcal G/\mathcal K$. Therefore $\Theta$ is the required complement to $\mathcal U/\mathcal K$. \medskip

Now assume that $\Theta \leq \mathcal G/\mathcal K$ is a complement to $\mathcal U/\mathcal K$ For $i=1,\ldots,d$ we can write $x_i\mathcal K$ uniquely as $x_i\mathcal K =z_i y_i$ with $z_i \in \Theta$ and $y_i \in \mathcal U/\mathcal K$. 
Define $\phi(x_i)=z_i$ ($i=1, \ldots, d$)  and extend $\phi$ to a homomorphism $\phi: \mathcal G \rightarrow \Theta \leq \mathcal G/\mathcal K$. By its definition $\phi$ satisfies P1. If $u \in \mathcal U$ then P1 with $y=u$ shows $\phi(u) \in \mathcal U/\mathcal K$. Thus $\phi(u) \in \phi(\mathcal G) \cap \frac{\mathcal U}{\mathcal K}= \Theta \cap \frac{\mathcal U}{\mathcal K} =\{1_{\frac{\mathcal H}{\mathcal K}}\}$ and hence $u \in \ker \phi$. Therefore P2 holds as well.
$\square$ 
\medskip

The maps $\phi: \mathcal G \rightarrow \frac{\mathcal G}{\mathcal K}$ which satisfy property P1 from Proposition \ref{split} are in bijective correspondence with maps  $h : \mathcal G \rightarrow S=\mathcal U /\mathcal K$  defined by $\phi(g)=h(g) g\mathcal K$ for each $g \in \mathcal G$. 
We observe that $\phi$ is a group homomorphism if and only if \[ h(g_1g_2)=h(g_1) \ ^{g_1} h(g_2) \quad \forall g_1,g_2 \in \mathcal G, \] that is, $h$ 
is a derivation from $\mathcal G$ into $S$ (where, as usual, $\mathcal G$ acts on $S$ via the projection  $\mathcal G \rightarrow \Gamma=\mathcal G/\mathcal U$). 

In addition $\phi$ satisfies property P2 if and only if $h(u)u\mathcal K  =1_{\mathcal U/\mathcal K}$ for all $u \in \mathcal U$, which in the additive group $(S,+)$ becomes  $h(u)=-f(u\Phi(\mathcal U))$ for each $u \in \mathcal U$.

In summary, the above discussion shows that $\mathcal U/\mathcal K$ is complemented in $\mathcal G/\mathcal K$ if and only if there is a derivation $h: \mathcal G \rightarrow S$ such that \begin{equation} \label{cond} h(u)= -f(u\Phi(\mathcal U) ), \quad \forall u \in \mathcal U. \end{equation}

A derivation $h: \mathcal G \rightarrow S$ is specified uniquely by its values $s_i:=h(x_i) \in S$ on the generating set $x_1, \ldots, x_{d}$ of $\mathcal G$ and then $h$ is determined by 
\[ h(w)= \sum_{i=1}^{d} \overline{\frac{\partial w}{\partial x_i}} s_i \quad \forall w \in \mathcal G. \]
Here $\frac{\partial w}{\partial x_i} \in \mathbb F_p \mathcal G$ are the Fox derivatives of the group algebra $\mathbb F_p \mathcal G$
and $y \mapsto \bar y \in \mathbb F_p \Gamma $ is the reduction $\mathbb F_p \mathcal G \rightarrow \mathbb F_p \Gamma$ of $\mathbb F_p \mathcal G$ modulo its ideal generated by $\mathcal U-1$.
For a proof see \cite{B}, \S IV.2, Exercise 3. 

Recall the description of $\rho: \mathcal U_p^{ab} \rightarrow (\mathbb F_p \Gamma)^{d}$ in (\ref{rho}):

\[ \rho(u\Phi(\mathcal U))= \left (\overline {\frac{\partial u}{\partial x_1}}, \ldots, \overline {\frac{\partial u}{\partial x_{d}}}\right ), \quad \forall u \in \mathcal U.\]

Let $\theta: (\mathbb F_p \Gamma)^{d} \rightarrow S$ be the homomorphism defined by $\theta(a_1, \ldots,a_{d})= \sum_{i=1}^{d} a_is_i$ for all $a_i \in \mathbb F_p\Gamma$. 

Combining the formulas for $h$ and $\rho$ above we obtain $h(u)=\theta \circ \rho(u \Phi(\mathcal U))$ for all $u \in \mathcal U$. Hence the condition (\ref{cond}) required from the derivation $h$ is equivalent to $-f(x)= \theta \circ \rho(x)$ for each $x \in \mathcal U_p^{ab}$ i.e. $f= (-\theta) \circ \rho$. Proposition \ref{lift} follows. $\square$

\end{document}